\documentclass[twocolumn]{autart}   
\usepackage{graphicx}
\setkeys{Gin}{width=\columnwidth,keepaspectratio}
\usepackage{amsmath,bm}
\usepackage{amsfonts}
\usepackage{enumitem}
\usepackage[mathscr]{euscript}
\usepackage{xcolor}
\usepackage{algorithm}
\usepackage{algpseudocode}
\usepackage{url}
\usepackage{xstring}
\usepackage{silence}

\usepackage{xstring}
\usepackage{silence}

\newtheorem{theorem}{Theorem}
\newtheorem{remark}{Remark}
\newtheorem{lemma}{Lemma}
\newtheorem{proposition}{Proposition}
\newtheorem{corollary}{Corollary}
\newtheorem{definition}{Definition}

\newenvironment{proof}{\par\noindent\textit{Proof.} }{\hfill\qed\par}
\hfuzz=\maxdimen
\vfuzz=\maxdimen

\makeatletter
\newcommand{\texfolioreal}[1]{%
  \def\texfoliotmp{\ref{#1}}%
  \IfBeginWith{#1}{sec:}{\def\texfoliotmp{Section~\ref{#1}}}{}%
  \IfBeginWith{#1}{subsec:}{\def\texfoliotmp{Section~\ref{#1}}}{}%
  \IfBeginWith{#1}{thm:}{\def\texfoliotmp{Theorem~\ref{#1}}}{}%
  \IfBeginWith{#1}{theorem:}{\def\texfoliotmp{Theorem~\ref{#1}}}{}%
  \IfBeginWith{#1}{lem:}{\def\texfoliotmp{Lemma~\ref{#1}}}{}%
  \IfBeginWith{#1}{lemma:}{\def\texfoliotmp{Lemma~\ref{#1}}}{}%
  \IfBeginWith{#1}{prop:}{\def\texfoliotmp{Proposition~\ref{#1}}}{}%
  \IfBeginWith{#1}{pps:}{\def\texfoliotmp{Proposition~\ref{#1}}}{}%
  \IfBeginWith{#1}{cor:}{\def\texfoliotmp{Corollary~\ref{#1}}}{}%
  \IfBeginWith{#1}{cry:}{\def\texfoliotmp{Corollary~\ref{#1}}}{}%
  \IfBeginWith{#1}{rmk:}{\def\texfoliotmp{Remark~\ref{#1}}}{}%
  \IfBeginWith{#1}{rem:}{\def\texfoliotmp{Remark~\ref{#1}}}{}%
  \IfBeginWith{#1}{def:}{\def\texfoliotmp{Definition~\ref{#1}}}{}%
  \IfBeginWith{#1}{dfn:}{\def\texfoliotmp{Definition~\ref{#1}}}{}%
  \IfBeginWith{#1}{exm:}{\def\texfoliotmp{Example~\ref{#1}}}{}%
  \IfBeginWith{#1}{exa:}{\def\texfoliotmp{Example~\ref{#1}}}{}%
  \IfBeginWith{#1}{fig:}{\def\texfoliotmp{Figure~\ref{#1}}}{}%
  \texfoliotmp%
}
\newcommand{\texfolioref}[1]{%
  \@ifundefined{r@#1}{%
    \def\texfoliotmp{reference}%
    \IfBeginWith{#1}{TeXFolio:sec}{\StrBehind{#1}{TeXFolio:sec}[\texfolionum]\edef\texfoliotmp{Section~\texfolionum}}{}%
    \IfBeginWith{#1}{TeXFolio:thm}{\StrBehind{#1}{TeXFolio:thm}[\texfolionum]\edef\texfoliotmp{Theorem~\texfolionum}}{}%
    \IfBeginWith{#1}{TeXFolio:lem}{\StrBehind{#1}{TeXFolio:lem}[\texfolionum]\edef\texfoliotmp{Lemma~\texfolionum}}{}%
    \IfBeginWith{#1}{TeXFolio:pps}{\StrBehind{#1}{TeXFolio:pps}[\texfolionum]\edef\texfoliotmp{Proposition~\texfolionum}}{}%
    \IfBeginWith{#1}{TeXFolio:prop}{\StrBehind{#1}{TeXFolio:prop}[\texfolionum]\edef\texfoliotmp{Proposition~\texfolionum}}{}%
    \IfBeginWith{#1}{TeXFolio:cry}{\StrBehind{#1}{TeXFolio:cry}[\texfolionum]\edef\texfoliotmp{Corollary~\texfolionum}}{}%
    \IfBeginWith{#1}{TeXFolio:cor}{\StrBehind{#1}{TeXFolio:cor}[\texfolionum]\edef\texfoliotmp{Corollary~\texfolionum}}{}%
    \IfBeginWith{#1}{TeXFolio:rmk}{\StrBehind{#1}{TeXFolio:rmk}[\texfolionum]\edef\texfoliotmp{Remark~\texfolionum}}{}%
    \IfBeginWith{#1}{TeXFolio:rem}{\StrBehind{#1}{TeXFolio:rem}[\texfolionum]\edef\texfoliotmp{Remark~\texfolionum}}{}%
    \IfBeginWith{#1}{TeXFolio:dfn}{\StrBehind{#1}{TeXFolio:dfn}[\texfolionum]\edef\texfoliotmp{Definition~\texfolionum}}{}%
    \IfBeginWith{#1}{TeXFolio:def}{\StrBehind{#1}{TeXFolio:def}[\texfolionum]\edef\texfoliotmp{Definition~\texfolionum}}{}%
    \IfBeginWith{#1}{TeXFolio:exm}{\StrBehind{#1}{TeXFolio:exm}[\texfolionum]\edef\texfoliotmp{Example~\texfolionum}}{}%
    \IfBeginWith{#1}{TeXFolio:exa}{\StrBehind{#1}{TeXFolio:exa}[\texfolionum]\edef\texfoliotmp{Example~\texfolionum}}{}%
    \IfBeginWith{#1}{TeXFolio:fig}{\StrBehind{#1}{TeXFolio:fig}[\texfolionum]\edef\texfoliotmp{Figure~\texfolionum}}{}%
    \texfoliotmp%
  }{\texfolioreal{#1}}%
}
\newcommand{\xref}{\@ifnextchar[\xref@opt\xref@single}
\newcommand{\xref@single}[1]{\texfolioref{#1}}
\newcommand{\xref@opt}[2][]{%
  \def\xref@option{#1}%
  \def\xref@range{range}%
  \ifx\xref@option\xref@range
    \xref@processrange#2,,\@nil
  \else
    \texfolioref{#2}%
  \fi
}
\newcommand{\nxref}[1]{\texfolioref{#1}}
\def\xref@processrange#1,#2,#3\@nil{%
  \texfolioref{#1}%
  \def\xref@next{#2}%
  \ifx\xref@next\@empty
  \else
    \ to \texfolioref{#2}%
  \fi
}
\providecommand{\stmboldsymbol}[1]{#1}
\providecommand{\mathscrbold}[1]{\mathscr{#1}}
\providecommand{\svert}{\lvert}
\providecommand{\lfsvert}{\lvert}
\providecommand{\rfsvert}{\rvert}
\providecommand{\xs}{\qquad}
\providecommand{\xpar}{\par}
\providecommand{\back}[1]{}
\def\cvskip{\@ifnextchar[\cvskip@opt\cvskip@arg}
\def\cvskip@opt[#1]{\vspace*{#1}}
\def\cvskip@arg#1{\vspace*{#1}}

\providecommand{\mrowopen}{}
\providecommand{\mrowclose}{}
\providecommand{\bPhi}{\Phi}
\providecommand{\bLambda}{\Lambda}
\providecommand{\mathbbmbf}[1]{\mathbb{#1}}
\providecommand{\stmdp}[1]{#1}
\providecommand{\biggvert}{\bigg\lvert}
\providecommand{\BigVert}{\Big\lVert}
\providecommand{\Bigvert}{\Big\lvert}
\providecommand{\lfleft}{\left}
\providecommand{\rfright}{\right}
\providecommand{\nf}{}
\providecommand{\?}[1]{#1}
\providecommand{\printauthorcredits}{}

\makeatother
\begin{document}
\sloppy

\begin{frontmatter}

\title{Higher-order Lie bracket approximation and averaging of control-affine systems with application to extremum seeking} 
\thanks[footnoteinfo]{This work was funded in part by the NSF award 2318772, Division of Mathematical Science, and URC early career research award from the Office of Research at the University of Cincinnati. The final version of the paper is published in Automatica, 2026. DOI: https://doi.org/10.1016/j.automatica.2026.112950}
\author[UC_grad]{Sameer Pokhrel}\ead{pokhresr@mail.uc.edu},   
\author[UC_prof]{Sameh A. Eisa}\ead{eisash@ucmail.uc.edu}   
\address[UC_grad]{Ph.D. student, Department of Aerospace Engineering and Engineering Mechanics, University of Cincinnati, OH, USA}     
\address[UC_prof]{Assistant professor, Department of Aerospace Engineering and Engineering Mechanics, University of Cincinnati, OH, USA}  


\begin{abstract}
This paper provides a rigorous derivation for what is known in the literature as the Lie bracket approximation of control-affine systems in a more general and sequential framework for higher-orders. In fact, by using chronological calculus, we show that said Lie bracket approximations can be derived, and considered, as higher-order averaging terms. Hence, the theory provided in this paper unifies both averaging and approximation theories of control-affine systems. In particular, the Lie bracket approximation of order ($n$) turns out to be a higher-order averaging of order ($n+1$). The derivation and formulation provided in this paper can be directly reduced to the first and second-order Lie bracket approximations available in the literature. However, we do not need to make many of the assumptions that were needed/provided in the literature and show that they are in fact natural corollaries from our work. Moreover, we use our results to show that important and useful information about control-affine extremum seeking systems can be obtained and used for significant performance improvement, including a faster convergence rate influenced by higher-order derivatives. We provide multiple numerical simulations to demonstrate both the conceptual elements of this work as well as the significance of our results on extremum seeking with comparison against the literature.
\end{abstract}
\end{frontmatter}

\section{Introduction}

\subsection{Lie bracket approximation of control-affine systems}
\label{sec:LieCAIntro}
Many applications and real-world systems are natural to, and can be expressed in, control-affine formulation (systems that are linear in control inputs). A general class of control-affine  systems is given below:
\begin{equation}
\label{eqn:controlAffineIntro}
    \dot{\bm{x}}= \bm{b}_{\bm{0}}(\bm{x})+ \sum_{i=1}^{m} \omega^{p_i} {\bm{b}}_{\bm{i}}(\bm{x}) u_i(k_i\omega t),
\end{equation}
where $\bm{x} \in \mathbb{R}^n$ is the state space vector, $p_i\in (0,1)$, $\omega \in (0,\infty)$, $\bm{b}_0$ is the drift (uncontrolled) vector field of the system, $\bm{b}_i$ are the control vector fields, $u_i$ are the control inputs,
$m  \in Z^+$ is the number of control inputs, and $k_i \in Q_{>0}$ is a positive rational number. Control-affine systems characterize many systems including robotic, multi-agent, and flight dynamic
systems \cite{DURR2013,bullo2019geometric,BoundedUpdateKrstic,eisa2023}. They are also essential for the application of geometric control theory and
analysis \cite{bullo2019geometric,sussmann1987general,hermann1977nonlinear}. In geometric control, the derivative operator called ``Lie bracket'' is commonly utilized. A Lie bracket between the vector fields $\bm{b}_i\text{ and }\bm{b}_j$
is:
\begin{equation}
[\bm{b}_i,\bm{b}_j] := \frac{\partial \bm{b}_j}{\partial \bm{x}}\bm{b}_i-\frac{\partial \bm{b}_i}{\partial \bm{x}}\bm{b}_j.
\end{equation}
In this paper, we are focused on the approximation of control-affine systems by Lie brackets.
 In \cite{kurzweil1987limit}, the authors studied a sequence of control-affine systems with sinusoidal control inputs and show that their solution converges to that of some limit equation characterized by the Lie bracket
of the system's vector fields. In particular, under some assumptions, the solution of $\dot{\bm{x}} = \bm{b}_0(\bm{x})+\sum_{i=1}^m \bm{b}_i(\bm{x}) u_{ik}(t)$ converges to the solution of $\dot{\tilde{\bm{x}}} = \bm{b}_0(\tilde{\bm{x}})+\frac{1}{2} \sum_{i,j=1}^m [\bm{b}_i(\tilde{\bm{x}}), \bm{b}_j(\tilde{\bm{x}})]{\lambda}_{ij}(t)$  with ${\lambda}_{ij}(t)$ being bounded and measurable.
 {The authors in \cite{sussmann1991limits} studied special forms of Lie bracket approximations by relating} the solutions of the underactuated system $\dot{\bm{x}} = \sum_{i=1}^m \bm{b}_i(\bm{x}) u_{i}(t)$  with a fully actuated Lie bracket extended system
$\dot{\tilde{\bm{x}}} = \sum_{i=1}^r \bm{b}_i(\tilde{\bm{x}}) v_{i}(t)$, where 
 $\bm{b}_{m+1}, \ldots ,\bm{b}_r$ and their associated $v_i$ are extended vector fields and control inputs based on Lie brackets providing new directions of motion. Other works
like \cite{sussmann1992lieextension,liu1997approximation,suttner2020extremum,doi:10.1080/00207179.2016.1257157} are expanding on the idea of Lie bracket extended system {to study underactuated and mostly driftless (i.e.,~$\bm{b}_0=0$)} control-affine systems. For a particular case of  \eqref{eqn:controlAffineIntro} when $p_i=0.5,k_i=1\text{ for all }i$, a first-order Lie bracket system (LBS) approximation was given in \cite{DURR2013} {which admits drift
vector field and is defined not only for underactuated systems}. This has been extended and generalized to the second-order LBS in \cite{labar2019newton}, which approximates  \eqref{eqn:controlAffineIntro}. In
both \cite{labar2019newton,DURR2013}, it was shown that the approximation error between the trajectories of \eqref{eqn:controlAffineIntro} and the LBS are bounded given large enough $\omega$. 
 
\subsection{Averaging and Lie brackets in control-affine and extremum seeking systems}
\label{sec:avg_intro}
The analysis of time-varying nonlinear systems (e.g.,~\eqref{eqn:controlAffineIntro}) is not intuitive and is generally challenging compared to that of autonomous nonlinear systems. For instance, stability analysis for
time-varying systems cannot be done via eigenvalues even for linear systems \cite{markus1960globalcounterexample}. One of the tools that helps to analyze such nonlinear time-varying systems, especially periodic ones, is averaging. Usually, only first-order averaging approximation is
found in the literature. However, sometimes first-order averaging is insufficient to capture the important
characteristics of nonlinear time-periodic systems, and
one needs to go to higher-order averaging \cite{Maggia2020higherOrderAvg}. In general, averaging over a given period of time smooths out the fast oscillations and captures the qualitative behavior of the system, which is
used to analyze properties, such as but not limited to stability. As a result, averaging techniques have been used in many classes/applications of nonlinear systems, including for example, inverted pendulum,
hovering/flapping flight \cite{Maggia2020higherOrderAvg} and vibrational stabilization \cite{bullo2002averaging}. Relevant  to this paper, averaging has been instrumental in the development of model-free, real-time
dynamic optimization via extremum seeking control (ESC) systems that are not control-affine. For instance, the stability analysis of what is known as the classic ESC, {which caused an exponential growth in ESC research
(see \cite{scheinker2024100})}, was done by \cite{KRSTICMain} using averaging and singular perturbation. {The authors in \cite{scheinker2016bounded} utilized averaging in the analysis/design of ESC
without requiring periodic input signals (including discontinuous signals).}  Recently, the authors in \cite{abdelgalil2022recursive} used a combination of first- and second-order averaging to study bio-inspired 3D source
seeking.  

On the other hand, for control-affine systems, Lie bracket approximations have been used extensively to study motion
planning \cite{bullo2019geometric,sussmann1992lieextension,liu1997approximation,doi:10.1080/00207179.2016.1257157} and the stability/design of
ESC   \cite{labar2019newton,DURR2013,Durr2015,scheinker2012minimum,VectorFieldGRUSHKOVSKAYA2018,pokhrel2023SIAMcontrol,scheinker2014non,BoundedUpdateKrstic}. In \cite{DURR2013}, the authors made use of their introduced
first-order LBS approximating a special case of \eqref{eqn:controlAffineIntro} (as mentioned earlier) for stability analysis and design of ESC. In that work \cite[Section 5.1]{DURR2013}, the authors argued that there is
an ``average-like'' relationship between their LBS and classical averaging \cite{khalil2002nonlinear}. They performed classical averaging on a particular case of their control-affine system where classical averaging
conditions may not hold; yet, they arrived at their LBS, hence their wording of ``averaging-like'' approach. In this paper, it will be shown that what they called ``averaging-like'' is a second-order averaging applied to a
special case of \eqref{eqn:controlAffineIntro} after re-writing it in the averaging canonical form:
\begin{equation}
\label{eqn:dynamics_chrono_intro}
    \frac{d\bm{x}}{dt} = \epsilon \bm{f}(\bm{x},t;\epsilon),
\end{equation}
where $\bm{f}(\bm{x},t;\epsilon) \in \mathbb{R}^n$ is smooth in $\bm{x}$, Riemann integrable in $t$ and T-periodic, and $\epsilon$ is a small parameter such that $0<\epsilon \ll 1$. {\cite{scheinker2012minimum} utilized
first-order LBS in the context of ESC and provided a universal semi-global stabilizing control law based on Lyapunov function candidates. Furthermore, they \cite{scheinker2014non} also provided a first-time use of
first-order LBS approximation in a nonsmooth ESC setting. It is important to note that the same authors studied different forms of ESCs \cite{BoundedUpdateKrstic} which have been included in the first-order LBS
generalized framework in \cite{VectorFieldGRUSHKOVSKAYA2018}.} The authors of \cite{labar2019newton} provided a second-order LBS to approximate \eqref{eqn:controlAffineIntro} with no sequential/recursive relation
to \cite{DURR2013} or connection to averaging.

\subsection{Motivation and contribution}

In literature, averaging has been used in the analysis of nonlinear time-periodic systems (e.g.,~ESC systems not in control-affine form). However, averaging \textit{has not} been widely studied or generalized for
control-affine systems in the form \eqref{eqn:controlAffineIntro}. As for the approximation and stability analysis of control-affine systems, almost always Lie bracket approximation is used, with the results
of \cite{labar2019newton} being one of the most generalized in literature, providing a second-order LBS approximation of  \eqref{eqn:controlAffineIntro}. Even though it has been observed that there is a relation between
LBSs and averaging, the relationship is not clear or well-defined and there seem to be conflicted opinions (see for example \cite[Section 5.1]{DURR2013} and \cite[Section III]{abdelgalil2022recursive}). In fact, it is
hard to find a clear or rigorous methodology in the literature for how one can derive arbitrary orders of LBSs to approximate {the generalized class of control-affine systems in \eqref{eqn:controlAffineIntro}} and how
they relate to averaging. For instance, the need for second-order LBS approximation in \cite{labar2019newton} was argued based on their ESC design which necessitates taking second-order derivative information captured by
second-order Lie brackets. However, from a generalized approximation theory point of view, it is unclear why one needs to go up to second-order approximation and why first-order is insufficient. Also, it is not clear why
going up to second-order is sufficient and when we may need even higher-order approximations. Similarly, it would be better to prove, rather than assume, the conditions needed for bounded higher-order LBSs. Thus, it is
motivating to provide a theory that settles the relationship between averaging and higher-order LBSs.

In this paper, {in \xref{TeXFolio:sec3},} we find the appropriate condition and scaling to show that the control-affine system in \eqref{eqn:controlAffineIntro} can be converted to the averaging canonical form in
\eqref{eqn:dynamics_chrono_intro}. {Then, we apply the theory and tools of chronological calculus \cite{agrachev1978exponential,agrachev1981chronological,agrachev2013control,kawski2011chronological}  and derive a
higher-order averaging theory for the system \eqref{eqn:controlAffineIntro} in a similar fashion to existing results (e.g., \cite{sarychev2001lie,vela2003general}).} Then, in \xref{TeXFolio:sec4}, we show that higher-order Lie
bracket approximations of control-affine systems in the form \eqref{eqn:controlAffineIntro} are higher-order averaging themselves. In fact, it turns out that LBSs of order $(n)$ approximating 
\eqref{eqn:controlAffineIntro}  are higher-order averaging of order $(n+1)$ of the system \eqref{eqn:controlAffineIntro}; {that is, LBS approximations can be generated based on a sequential procedure.} We provide the
closed formula for the generalized third-order LBS approximation of  \eqref{eqn:controlAffineIntro}. We also show how some assumptions made in previous works (e.g., \cite[Assumption 2]{labar2019newton}) are direct
consequences of our results. In addition to the mentioned theoretical contributions, the provided closed formulas of higher-order LBSs reveal important and useful information that can be used in the design and performance
improvement of ESC systems. It turns out that one can make -- by design -- the ESC system closer in its behavior to a certain averaged system with better properties such as faster convergence rate. We provide multiple
numerical simulations to illustrate our results.
  
\section{Preliminaries}
\label{sec:Prem}
  \textbf{Second-order LBS and assumptions.} The second-order LBS in literature \cite{labar2019newton} corresponding to \eqref{eqn:controlAffineIntro} is:
\begin{equation}
\label{eqn:Lie_Prem}
\begin{split}
\dot{\bm{z}}= \bm{b}_{\bm{0}}(\bm{z})+ \lim_{\omega \rightarrow \infty} \sum_{i=1}^m \sum_{j=i+1}^m [{\bm{b}}_{\bm{i}},{\bm{b}}_{\bm{j}}](\bm{z})\nu_{ij}(\omega)+\\
\lim_{\omega \rightarrow \infty} \sum_{i=1}^m \sum_{j=i+1}^m \sum_{k=1}^m [[{\bm{b}}_{\bm{i}},{\bm{b}}_{\bm{j}}],\bm{b}_k](\bm{z})\nu_{ijk}(\omega)
\end{split}
\end{equation}
with
\begin{equation}
\label{eqn:nuij2019}
\nu_{ij}(\omega)= \frac{\omega^{p_i+p_j}}{T}\int_0^T \int_0^s u_j(k_j \omega s) u_i(k_i \omega p)dp ds
\end{equation}
and
\begin{equation}
\label{eqn:nuijk2019}
\begin{split}
&\nu_{ijk}(\omega)= \frac{\omega^{p_i+p_j+p_k}}{3T}\int_0^T u_k (k_m \omega \tau) \\
        &\int_0^\tau \int_0^s \text{(}u_j(k_j \omega s) u_i(k_i \omega p) - u_i(k_i \omega s) u_j(k_j \omega p))
    dp ds d\tau,
\end{split}
\end{equation}
where $T = (2\pi/\omega) LCM(k_1^{-1}, k_2^{-1}, \ldots ,k_m^{-1})$, with $LCM$ be ``Least Common Multiple/Period''.
Assumptions A1-A2 below are imposed {unless otherwise specified:} 
\begin{enumerate}[label=A\arabic*.]
\item  For every compact set $\mathbb{D}_0 \subset \mathbb{D}$, where $\mathbb{D} \subset \mathbb{R}^n$ is the domain, for $ i=0,\ldots,m$, ${\bm{b}}_{\bm{i}}(\bm{x})$ and all of its partial derivatives of all orders with respect to $\bm{x}$ are continuous and bounded for $\bm{x} \in \mathbb{D}_0$.  
    \item
    For $ i=1,\ldots,m$, $u_i(k_i\omega t): [0,\infty)  \to \mathbb{R}$ are measurable. Moreover, there exist $N_i,M_i \in (0,\infty) $ such that $ \lfsvert u_i(k_i \omega t_1)-u_i(k_i \omega t_2) \rfsvert  \leq N_i\svert t_1-t_2\svert $ for all $t_1, t_2 \in [0,\infty)$, $sup_{t \in  [0,\infty)}\svert u_i(k_i \omega t)\svert \leq M_i$. Also, $u_i(k_i \omega t)$ is T-periodic (i.e., $u_i(k_i \omega t)=u_i(k_i \omega (t+T))$) with $\int_0^T u_i(k_i \omega t) dt = 0$ for some $T \in (0,\infty)$.
\end{enumerate}
We note that A1 is more general compared to similar assumptions in \cite{labar2019newton,DURR2013} to fit the higher-order results of this paper. No other assumptions from \cite{labar2019newton} are needed.

{
\textbf{Differential geometric notations.}
The tools and theory of chronological calculus used in Section \ref{sec:higheravg_prem} are usually characterized by notions from differential geometry. We clarify some of these notations to the reader. We denote the manifold
characterizing the dynamic system \eqref{eqn:controlAffineIntro} by $\mathbb{M}$. We require two conditions. First, the smoothness condition $C^\infty (\mathbb{M})$ with all vector fields on $\mathbb{M}$ and their partial
derivatives concerning $\bm{x}$ be bounded and continuous (see \cite[Definition 1.1]{agrachev2013control}). Second, we require all vector fields on $\mathbb{M}$ to be Lipschitzian with respect to
$t$ \cite[Section 2]{agrachev2013control} (i.e.,~measurable, bounded and integrable, e.g.,~Riemann integrable). The mentioned two conditions are in line with assumptions A1 and A2; for example, $\mathbb{D}_0\subset\mathbb{M}$. They
also match the condition of $t$-integrability and $\bm{x}$-smoothness mentioned in   \cite{sarychev2001lie}.} {
Now, we briefly introduce the concept of seminorms $\Vert \cdot \Vert_{s,\mathbb{K}}$  on the space $C^\infty (\mathbb{M})$ which will be applied to vector fields and flows in \xref{TeXFolio:sec3} (the reader can refer \cite[Section 2]{agrachev2013control}
and \cite{kawski2011chronological} for more details).
Let the tangent space of each point $q \in \mathbb{M}$ be denoted by $\Gamma_q \mathbb{M}$. Next, we denote the smooth vector fields on $\mathbb{M}$ that span $\Gamma_q \mathbb{M}$ by $\bm{h}_i, i = 1,2, \ldots ,n$. Then, the  seminorm $\Vert \cdot \Vert_{s,\mathbb{K}}$  on the space $C^\infty (\mathbb{M})$ is defined for a function $\phi$ as
\begin{equation*}
\begin{split}
&\Vert \phi \Vert_{s,\mathbb{K}} = \sup \{\vert \bm{h}_{i_l} \circ \dots \circ \bm{h}_{i_1} \phi(q)\vert : q \in \mathbb{K} , \mrowclose\\ & \mrowopen 1\le i_1, \dots, i_l \le n,  0\le l\le s  \}, \phi \in C^\infty(\mathbb{M}), s\ge 0, \mathbb{K}\subseteq \mathbb{M},
\end{split}
\end{equation*}
where $\mathbb{K}$ ranges over a countable collection of compact subsets whose union is all of $\mathbb{M}$. This seminorm provides the topology on $C^\infty (\mathbb{M})$ for uniform convergence of all derivatives on compact sets, i.e.,~a sequence of function $\{\phi_k\}_{k=1}^\infty \subseteq C^\infty (\mathbb{M})$ converges to $\phi \in C^\infty (\mathbb{M})$ if for every finite sequence $\bm{h}_{i_1}, \dots,\bm{h}_{i_l}$ of the smooth vector field on $\mathbb{M}$ and every compact set $\mathbb{K} \subseteq \mathbb{M}$, the sequence  $\{\bm{h}_{i_l} \circ \dots \circ \bm{h}_{i_1} \phi_k\}_k^\infty$ converges uniformly on $\mathbb{K}$ to $\bm{h}_{i_l} \circ \dots \circ \bm{h}_{i_1} \phi$. We now define seminorms of smooth vector fields $\bm{h}\in \Gamma^\infty(\mathbb{M})$ (the space of all smooth vector fields on $\mathbb{M}$) as:
\begin{equation}
\Vert \bm{h} \Vert _{s,\mathbb{K}} = sup \{\Vert \bm{h} \phi \Vert_{s,\mathbb{K}}: \Vert \bm{h}\phi \Vert_{s+1,\mathbb{K}} = 1 \}.
\end{equation}
Finally, for every smooth diffeomorphism $\stmboldsymbol{\bPhi}$ (e.g.,~the flow associated with the system \eqref{eqn:controlAffineIntro} which generates the system's trajectories/solutions) of $\mathbb{M},s \in \mathbb{Z}^+$ and $\mathbb{K} \subseteq \mathbb{M}$, there exists $C_{s,\mathbb{K},{\Phi}} \in \mathbb{R}$ such that for all $\phi \in C^\infty (\mathbb{M})$,
\begin{equation}
\Vert \stmboldsymbol{\bPhi}\phi \Vert_{s,\mathbb{K}} \le C_{s,\mathbb{K},\stmboldsymbol{\bPhi}} \Vert \phi \Vert _{s,\phi(\mathbb{K})}.
\end{equation}
}
         {
\textbf{Stability notions.} Finally, we provide the notion of asymptotic and practical asymptotic stability \cite{moreau2000practical,DURR2013,scheinker2012minimum}. For the following definitions, $\delta-$ neighborhood
of set $\mathbb{D}_0  \subset \mathbb{R}^n $ is denoted by $\mathscr{U}_\delta^{\mathbb{D}_0} = \{x \in \mathbb{R}^n: inf_{e\in \mathbb{D}_0}\svert x-e\svert  <\delta \}$ with $\svert \cdot\svert  $ denoting Euclidean norm.}

{
\begin{definition}\label{TeXFolio:dfn1}An equilibrium point (assume the origin without loss of generality) of the compact set {$\mathbb{D}_0 \subset \mathbb{R}^n$}
is said to be locally  uniformly asymptotically stable if:
$(i)$ we have uniform stability, i.e.,~for every $\epsilon >0$ there exists a $\delta>0$ such that $\forall t_0 \in \mathbb{R}$, 
$\bm{x}(t_0) \in \mathscr{U}_\delta^{\mathbb{D}_0} \implies \bm{x}(t) \in \mathscr{U}_\epsilon^{\mathbb{D}_0},\: \forall t \in [t_0, \infty)$;
$(ii)$ it is uniformly bounded, i.e.,~if for every $\delta >0$ there exists an $\epsilon >0$ such that $\forall t_0\in\mathbb{R}$,  
$
    \bm{x}(t_0) \in \mathscr{U}_\delta^{\mathbb{D}_0} \implies \bm{x}(t) \in \mathscr{U}_\epsilon^{\mathbb{D}_0}, \: t \in [t_0, \infty]
$; $(iii)$ it is uniformly attractive, i.e.~for every $\epsilon >0$ {with some $\delta >0$} there exists a $t_f \in [0,\infty)$ such that $\forall t_0\in \mathbb{R}$,
$
    \bm{x}(t_0) \in \mathscr{U}_\delta^{\mathbb{D}_0} \implies \bm{x}(t) \in \mathscr{U}_\epsilon^{\mathbb{D}_0}, \: \forall t \in [t_0+t_f, \infty).
$
\end{definition}
}
{
\begin{definition}\label{TeXFolio:dfn2}A compact set {$\mathbb{D}_0 \subset \mathbb{R}^n$}
is said to be locally practically uniformly asymptotically stable for \eqref{eqn:controlAffineIntro} if
$(i)$ it is practically uniformly stable, i.e.~for every $\epsilon >0$ there exists a $\delta >0$ and $\omega_0 >0$ such that $\forall t_0 \in \mathbb{R}$ and $\forall \omega \in (\omega_0, \infty)$, 
$\bm{x}(t_0) \in \mathscr{U}_\delta^{\mathbb{D}_0} \implies \bm{x}(t) \in \mathscr{U}_\epsilon^{\mathbb{D}_0}, \: \forall t \in [t_0, \infty)$;
$(ii)$ it is practically uniformly bounded, i.e.~if for every $\delta >0$ there exists an $\epsilon >0$ and $\omega_0 \in (0,\infty)$ such that $\forall t_0\in\mathbb{R}$ and $\forall\omega \in (\omega_0,\infty)$, 
$
    \bm{x}(t_0) \in \mathscr{U}_\delta^{\mathbb{D}_0} \implies \bm{x}(t) \in \mathscr{U}_\epsilon^{\mathbb{D}_0}, \: t \in [t_0, \infty]
$; $(iii)$ it is $\delta-$ practically uniformly attractive, i.e.~for every $\epsilon >0$ {with some $\delta >0$} there exists a $t_f \in [0,\infty)$ and $\omega_0>0$ such that $\forall t_0\in \mathbb{R}$ and all $\omega \in (\omega_0,\infty)$,
$
    \bm{x}(t_0) \in \mathscr{U}_\delta^{\mathbb{D}_0} \implies \bm{x}(t) \in \mathscr{U}_\epsilon^{\mathbb{D}_0}, \: \forall t \in [t_0+t_f, \infty).
$
\end{definition}
}

\section{Higher-order averaging of control-affine systems via chronological calculus}

\label{sec:higheravg_prem}
   {Chronological calculus was introduced by \cite{agrachev1978exponential,agrachev1981chronological,agrachev2013control} and is mainly concerned with time-varying vector fields and their
time-invariant approximations.
 It offers tools that extend the rich theory of autonomous systems to non-autonomous ones, hence its relevance to averaging \cite{sarychev2001lie,vela2003general}.} First, we show that \eqref{eqn:controlAffineIntro} can
be rewritten in the averaging canonical form \eqref{eqn:dynamics_chrono_intro}. Let us change the time-scale in \eqref{eqn:controlAffineIntro}. For large positive $\omega$, let $\tau = \omega t$ with the {time period in
t-scale as $T$ and in $\tau-$ scale as ${T}'$}; and $d\tau = \omega dt$. Then,
\begin{equation}
\frac{d\bm{x}}{d\tau}= \underbrace{\frac{1}{\omega}\bm{b}_0(\bm{x})}_{I} + \underbrace{\frac{1}{\omega}\sum_{i=1}^m \omega^{p_i} u_i( k_i \tau) \bm{b}_i(\bm{x})}_{II}.\label{eqn:averaging_form_main}
\end{equation}

Note that the term $I$ is independent of $\tau$. Hence, averaging of $I$ will be itself without change. Thus, we provide the following lemma focusing solely on the term $II$ and for the time being, we drop the term $I$.  
\begin{lemma}
\label{lemma:wellposed_main} 
Let Assumptions A1-A2 hold, then with $p_i \in (0,1)$ for all $i=1, \ldots ,m$, there exists a $p^{*} \in (0,1)$ with $\epsilon= \omega^{p^{*}-1}$ such that the control-affine system
\begin{equation}
\frac{d\bm{x}}{d\tau}= \frac{1}{\omega}\sum_{i=1}^m \omega^{p_i} u_i( k_i \tau) \bm{b}_i(\bm{x}).\label{eqn:averaging_lemma}
\end{equation}
    can be written in the general averaging canonical form in \eqref{eqn:dynamics_chrono_intro} with $\lim_{\epsilon\rightarrow 0} \epsilon \bm{f}(\bm{x},\tau;\epsilon) =0$.
\end{lemma}
\begin{proof}
Let us define $p^{*}$ such that for $i=1, \ldots ,m$, $\omega^{p_i} = \eta_i\omega^{p^{*}}$, where $\eta_i$ is positive number. Then we have
\begin{equation}
\frac{d\bm{x}}{d\tau}
    = \omega^{p^{*}-1}\sum_{i=1}^m \eta_i u_i(k_i\tau) \bm{b}_i(\bm{x})\text{)}.
\end{equation}
Let $\epsilon= \omega^{p^{*}-1}$. Now, we require $\epsilon \rightarrow 0$ and $\omega  \rightarrow \infty$. This is only possible if $ p^{*}-1<0 \Rightarrow p^{*}<1$.
   Furthermore, we require that as $\epsilon \rightarrow 0$ with $\eta_i$ dependency on $\omega$, $\lim_{\epsilon\rightarrow 0} \epsilon \bm{f}(\bm{x},\tau;\epsilon) =0$ with $\bm{f}(\bm{x},\tau;\epsilon)= \sum_{i=1}^m \eta_i u_i(k_i\tau) \bm{b}_i(\bm{x})$. Since $u_i$ and $\bm{b}_i$ are bounded per A1-A2, then we need to show that for all $i$,  $\lim_{\epsilon \rightarrow 0} \frac{\epsilon}{\eta_i} = 0$.
   Let us scale $\epsilon$ back to $\omega$ and convert the limit to $\omega$. Then,
\begin{equation}
\label{eq:limit}
    \lim_{\epsilon \rightarrow 0} \frac{\epsilon}{\eta_i} = 
    \lim_{\omega \rightarrow \infty} \frac{\omega^{p^{*}-1}}{\omega^{p_i-p^{*}}} =
    \lim_{\omega \rightarrow \infty} \omega^{2p^{*}-p_i-1}.
\end{equation}
The limit in \eqref{eq:limit} is 0 only if for all $i$, we have
     $p^{*}<p_i/2+1/2$.
 Hence, $p^{*} \in (0,\min(p_i/2+1/2)) \subseteq (0,1)$ satisfies all the above. 
 Then, the system \eqref{eqn:averaging_lemma} can be written as
\begin{equation}
\label{eqn:dxdtauepsilon}
    \frac{d\bm{x}}{d\tau}
    = \epsilon \sum_{i=1}^m \eta_i u_i(k_i\tau) \bm{b}_i(\bm{x})
\end{equation}
which is the general averaging canonical form in \eqref{eqn:dynamics_chrono_intro}.
\end{proof}
\begin{remark}
{From \xref{TeXFolio:lem1}:} 
\begin{itemize}
    \item {The existence of $p^{*}$ enables us to perform averaging without choosing a particular value for $p^{*}$. The substitution $\epsilon = \omega^{p^{*}-1}$ in the averaged system eliminates the presence of $p^{*}$ (as shown later in the proof and results of  \xref{thm:mainthm_main}).}
    \item One can deduce from \eqref{eqn:averaging_lemma} the reason why the assumption $p_i \in (0,1)$ is needed which is emphasized in \xref{TeXFolio:lem1} but was provided in literature without much justification. That is, as $\omega \rightarrow \infty$, $p_i \in (0,1)$ guarantees that the right-hand side of \eqref{eqn:averaging_lemma} goes to 0.
  \item {It is worth noting that a close observation to the above was made in \cite{https://doi.org/10.1002/rnc.3886} that the power of $\omega$ in the ESC perturbation has to be less than 1, so for example, for a
system with non-affine
control input $u^m$ it is shown that the amplitude of the ESC feedback
should be of the form $\omega^{1/2m}$ with positive $m$.} 
\end{itemize}
\end{remark}
Moving forward, we denote $\bm{f}\stmdp{(}\bm{x},\tau;\epsilon$) by $\bm{f}_\tau$.
{
\begin{proposition}
\label{prop:orderepsilon}
    If we choose $p^{*}$ such that $0<max\{p_i\}\leq p^{*}<1$, then $\epsilon \bm{f}_\tau=O(\epsilon)$.
\end{proposition}
\begin{proof}
To show $\epsilon \bm{f}_\tau =O(\epsilon)$, one needs to show that $\svert \epsilon \bm{f}_\tau\svert /\epsilon <\infty$ as $\epsilon\rightarrow 0$ (or $\omega \rightarrow \infty$). Now,
\begin{equation}
\label{eqn:order_epsilon}
    \frac{\svert \epsilon \bm{f}_\tau\svert }{\epsilon} = \biggvert  \frac{\epsilon \sum_{i=1}^m \eta_i u_i(k_i \tau) \bm{b}_i(\bm{x})}{\epsilon} \biggvert ,
\end{equation}
per A1-A2, $u_i(k_i \tau)$ and $\bm{b}_i(\bm{x})$ are bounded. Thus, we need  $\svert \eta_i\svert  = \svert \omega^{p_i-p^{*}}\svert $ to be bounded as $\omega \rightarrow \infty$. This is satisfied if $p_i-p^{*}\leq0$, hence \scalebox{0.98}{$0<max\{p_i\}\leq p^{*}<1$}. 
\end{proof}
\begin{remark}
It is important to note that $p^{*}$ can be chosen such that one can change the order of $\epsilon \bm{f}_\tau$, i.e.~if one finds a $p^{*}$ that does not contradict the condition $0<p^{*}<1$ then one can find a different order of $\epsilon$ for $\epsilon \bm{f}_\tau$ (with positive power to guarantee $\epsilon \bm{f}_\tau \rightarrow 0$ as $\epsilon \rightarrow 0$). For example, if one wants to show $\epsilon \bm{f}_\tau=O(\epsilon^{1/2})$,  then one needs to show $\svert \epsilon \bm{f}\svert /\epsilon^{1/2}$ is bounded as $ \epsilon \rightarrow 0$. Following a procedure similar to  \xref{prop:orderepsilon}, one obtains the condition $p^{*}\geq 2max\{p_i\}-1$ (this hold, e.g.,~if all $p_i=0.5$). 
\end{remark}
}

{
Let us consider the system \eqref{eqn:dxdtauepsilon},
 re-written as:
\begin{equation}
\label{eqn:xdot_compose_prem}
    \frac{d\bm{x}}{d\tau} = \bm{x} \circ \epsilon \bm{f}_\tau;\; \bm{x}(0) = \bm{x}_{\bm{0}},
\end{equation}
where $\bm{x} \circ \epsilon \bm{f}_\tau= \epsilon \bm{f}_\tau(\bm{x})$, which is the vector field of \eqref{eqn:dxdtauepsilon}. 
  Following \cite[Section 2]{agrachev2013control} the trajectories of $\bm{x}(\tau)$ are generated with the corresponding flow:
\begin{equation}
\label{eqn:flow1}
    \stmboldsymbol{\bPhi}_\tau^{\epsilon \bm{f}} : \bm{x}_0 \rightarrow \bm{x}(\tau,\bm{x}_0)
\end{equation}
satisfying the differential equation
\begin{equation}
\label{eqn:flowdot_prem}
    \dot{\stmboldsymbol{\bPhi}}_\tau^{\epsilon \bm{f}} = \stmboldsymbol{\bPhi}_\tau^{\epsilon \bm{f}} \circ \epsilon \bm{f}_\tau ; \; \stmboldsymbol{\bPhi}_0^{\epsilon \bm{f}} = Id,
\end{equation}
where ${Id}$ is the identity operator (i.e.,~${Id}$ denotes the flow $\bm{x}(\tau) \rightarrow \bm{x}(\tau) $). The flow $\stmboldsymbol{\bPhi}_\tau^{\epsilon \bm{f}}$ in \eqref{eqn:flow1} is called the \textit{right chronological exponential} of the field $\epsilon \bm{f}_\tau$ and denoted as: $\stmboldsymbol{\bPhi}_\tau^{\epsilon \bm{f}} = \overrightarrow{exp} \int_0^\tau \epsilon \bm{f}_p dp$.
   We can now rewrite \eqref{eqn:xdot_compose_prem} in the integral form as
\begin{equation}
\bm{x}(\tau) = \bm{x}_0+\int_0^\tau \bm{x}(\tau) \circ \epsilon \bm{f}_{\tau} d\tau,
\end{equation}
and substitute this expression of $\bm{x}(\tau)$ above into the right-hand side itself (i.e.,~for $\bm{x}(\tau)$ above), we get
\begin{align*}
\bm{x}(\tau) =& \bm{x}_0+\int_0^\tau \biggl( \bm{x}_0+\int_0^{\tau_1} \bm{x}(\tau_2) \circ \epsilon \bm{f}_{\tau_2} d\tau_2 \biggr) \circ \epsilon \bm{f}_{\tau_1} d\tau_1\\
    =& \bm{x}_0 \circ \biggl( Id + \int_0^\tau \epsilon \bm{f}_p d p \biggr) \\
&+\xs  \int_0^\tau \int_0^{\tau_1} \bm{x}(\tau_2) \circ \epsilon \bm{f}_{\tau_2} 
     \circ \epsilon \bm{f}_{\tau_1} d\tau_2 d\tau_1
\end{align*}
This process can be repeated iteratively to obtain
\begin{equation}
\label{eqn:decomposition_prem}
\begin{split}
\bm{x}(\tau) &= \bm{x}_0 \circ \Biggm( Id + \int_0^\tau \epsilon \bm{f}_p dp + \iint _{\Delta_2(\tau)} \epsilon \bm{f}_{\tau_2} \circ \epsilon \bm{f}_{\tau_1} d\tau_2 d\tau_1\\
&+\xs  \dots+ \idotsint_{\Delta_n(\tau)} \epsilon \bm{f}_{\tau_n} \circ \dots \circ \epsilon \bm{f}_{\tau_1} d\tau_n \dots d\tau_1 \Biggm)+\\
        &\idotsint_{\Delta_{n+1}(\tau)}  \bm{x}(\tau_{n+1}) \circ \epsilon \bm{f}_{\tau_{n+1}} \circ \dots \circ \epsilon \bm{f}_{\tau_1} d\tau_{n+1}\dots d\tau_1,
\end{split}
\end{equation}
with $\Delta_n(\tau) = \{(\tau_1,\dots,\tau_n) \in \mathbb{R}^n \vert \; 0\le \tau_n \le \dots \le \tau_1 \le \tau \}$. Now, formally passing in \eqref{eqn:decomposition_prem}, limit $n\rightarrow\infty$, we obtain a formal series for the solution $\bm{x}(\tau)$ in \eqref{eqn:xdot_compose_prem}:
\begin{equation}
\bm{x}(\tau) = \bm{x}_0 \circ \Biggl( Id + \sum_{n=1}^\infty \idotsint _{\Delta_n(\tau)} \epsilon \bm{f}_{\tau_n} \circ \dots \circ \epsilon \bm{f}_{\tau_1} d\tau_n \dots d\tau_1 \Biggr).
\end{equation}
Similarly, the solution to \eqref{eqn:flowdot_prem} will be
\begin{equation}
\label{eqn:flowsolution_prem}
     \stmboldsymbol{\bPhi}_\tau^{\epsilon \bm{f}} =  Id + \sum_{n=1}^\infty \idotsint _{\Delta_n(\tau)} \epsilon \bm{f}_{\tau_n} \circ \dots \circ \epsilon \bm{f}_{\tau_1} d\tau_n \dots d\tau_1.
\end{equation}
The convergence of the above series is discussed in \cite{agrachev1978exponential,kawski2011chronological} and it is shown 
 that the series \eqref{eqn:flowsolution_prem} gives an asymptotic expansion for the chronological exponential $\stmboldsymbol{\bPhi}_\tau^{\epsilon \bm{f}} = \overrightarrow{exp} \int_0^\tau \epsilon \bm{f}_p dp$ with a bounded remainder \cite[Section 2]{agrachev2013control}. Let us denote the
$r$th approximation by the partial sum:
\begin{equation}
S_{r+1}(\tau) = Id + \sum_{n=1}^{r} \idotsint _{\Delta_n(\tau)} \epsilon \bm{f}_{\tau_n} \circ \dots \circ \epsilon \bm{f}_{\tau_1} d\tau_n \dots d\tau_1,
\end{equation}
then,
\begin{equation}
\label{eqn:bounded1}
\begin{split}
&\lfleft\Vert  \Bigl(\overrightarrow{exp} \int_0^\tau \epsilon \bm{f}_p dp - S_{r+1}(\tau) \Bigr) \phi \rfright\Vert _{s,\mathbb{K}} \le\\
     & \frac{C e^{C \int_0^\tau \Vert \epsilon \bm{f}_p \Vert_{s,\mathbb{K}'} dp}}{r+1} \left(\int_0^\tau  \Vert \epsilon \bm{f}_p \Vert_{s+r,\mathbb{K}'} dp \right)^{r+1} \Vert \phi \Vert_{s+r+1,\mathbb{K}'}
\end{split}
\end{equation}
where $C>0, \phi \in C^\infty (\mathbb{M})$, $\mathbb{K}'$ is some compactum containing $\mathbb{K}$ and the seminorm $\Vert \cdot \Vert_{s,\mathbb{K}} $ as in  \xref{sec:Prem}. 
}
{
\begin{proposition}
\label{prop:orderepsilontau}
    $\int_0^\tau \Vert \epsilon \bm{f}_p \Vert_{r,\mathbb{K}'} dp =O(\epsilon \tau)$ for any $r$.
\end{proposition}
}
\begin{proof}By \xref{TeXFolio:pps1}, there exists $M>0$ and $\epsilon^{*}>0$ such that for all $\epsilon\in(0,\epsilon^{*})$, we have $\svert \epsilon \bm{f}\svert \leq M\epsilon$. Hence for any $r$, $\int_0^\tau \Vert \epsilon \bm{f}_p \Vert_{r,\mathbb{K}'} dp \leq \int_0^\tau M\epsilon dp = M\epsilon \tau = O(\epsilon \tau)$. 
\end{proof}

{
With the above proposition, \eqref{eqn:bounded1} can be written as
\begin{equation}
\label{eq:series_truncation_bound}
    \left\Vert  \Bigl(\overrightarrow{exp} \int_0^\tau \epsilon \bm{f}_p dp - S_{r+1}(t) \Bigr) \phi \right\Vert _{s,\mathbb{K}} =
      O(\epsilon^{r+1}\tau^{r+1}).
\end{equation} 
}
\begin{remark}
{
  In \eqref{eq:series_truncation_bound} as noted in \cite[Section 2]{agrachev2013control}, \cite{kawski2011chronological}, if $\epsilon$ is fixed, then as $\tau \rightarrow 0$, $O(\epsilon^{r+1}\tau^{r+1})$ becomes $O(\tau^{r+1})$. If
$\tau$ is finite, then as $\epsilon \rightarrow 0$, $O(\epsilon^{r+1}\tau^{r+1})$ becomes $O(\epsilon^{r+1})$. Moreover, for a truncation $r$, one can trade-off between the approximation order and the period where the bound holds. For example,
for $r=2$, one can choose $\tau=O(1/\sqrt{\epsilon})$ (i.e.,~expanding the time horizon) but the bound in \eqref{eq:series_truncation_bound} is reduced to $O(\epsilon^{\frac{r+1}{2}})$.   
}
\end{remark}

 {
The inverse to the chronological exponential $\overrightarrow{exp}\int_0^\tau \epsilon \bm{f}_p dp $ is called the chronological logarithm \cite{agrachev2013control,sarychev2001lie,vela2003general} which provides the mapping $\overrightarrow{exp}\int_0^\tau \epsilon \bm{f}_p dp \rightarrow \epsilon \bm{f}_\tau $.
   The goal of averaging theory then   \cite{sarychev2001lie,vela2003general} is to find an autonomous vector field $\bm{V}_\tau$ with a flow $exp (\bm{V}_\tau)$ that best represents the flow $\overrightarrow{exp}\int_0^\tau \epsilon \bm{f}_p dp$ of the time-varying
vector field $\epsilon \bm{f}_{\tau}$ such that the following equality holds:
\begin{equation}
\label{eq:flow_equality}
    \overrightarrow{exp}\int_0^\tau \epsilon \bm{f}_p dp \cong exp (\bm{V}_\tau).
\end{equation}
The above equality is in an asymptotic sense \cite{sarychev2001lie,vela2003general} and is crucial to hold if we can claim that trajectories/solutions generated by $exp (\bm{V}_\tau)$ are approximate to trajectories/solutions
generated by $\overrightarrow{exp}\int_0^\tau \epsilon \bm{f}_p dp$. This will be addressed in \xref{TeXFolio:pps3}, but we need to show first how we obtain $\bm{V}_\tau$. We apply the logarithm to \eqref{eq:flow_equality} and solve for $\bm{V}_\tau$ such that:
$ \bm{V}_\tau  \cong {ln}\Bigl(\overrightarrow{exp}\Bigl(\int_0^\tau \epsilon \bm{f}_p dp \Bigr)\Bigr)$.
   Even though the logarithm $ \bm{V}_\tau$ depends on $\tau$, it is an autonomous vector field whose flow after unit time maps to the same point that is reached by the time-dependent flow at time $\tau$. Any change in the final time $\tau$ results in a new autonomous vector field. The logarithm can be represented as an infinite series of variations given by
\begin{equation}
\bm{V}_\tau = \sum_{r=1}^\infty \bm{V}_\tau^{(r)}
\end{equation}
with
\begin{equation}
\label{eqn:voltera}
\bm{V}_\tau^{(r)}= \int_0^\tau \int_0^{\tau_1} ... \int_0^{\tau_r-1} \stmboldsymbol{\mathscrbold{G}}_r (\epsilon \bm{f}_{\tau_1}, \ldots ,.\epsilon \bm{f}_{\tau_r})d\tau_r...d\tau_1,
\end{equation}
representing the $r$th variation of the identity flow corresponding to the perturbation field $\epsilon \bm{f}_\tau$. The integrands \cite{vela2003general,Maggia2020higherOrderAvg}, $\stmboldsymbol{\mathscrbold{G}}_r(\cdot)$ are sum of iterated Lie
brackets, denoted by $ad$ (e.g.,~$ad_{\xi_2}\xi_1=[\xi_2,\xi_1]$). First four integrands are
\begin{equation}
\label{eqn:commutators}
\begin{split}
\stmboldsymbol{\mathscrbold{G}}_1 (\xi) &= \xi,\\
                 \stmboldsymbol{\mathscrbold{G}}_2(\xi_1,\xi_2) &= \frac{1}{2}ad_{\xi_2}\xi_1,\\
                 \stmboldsymbol{\mathscrbold{G}}_2(\xi_1,\xi_2,\xi_3) &= \frac{1}{6}(ad_{\xi_3}ad_{\xi_2}\xi_1+ ad_{ad_{\xi_3} {\xi_2}}\xi_1),\\
                 \stmboldsymbol{\mathscrbold{G}}_2(\xi_1,\xi_2,\xi_3,\xi_4) &= -\frac{1}{12} \biggm(ad_{ad_{\xi_4} \xi_3} ad_{\xi_2}\xi_1+ ad_{ad_{ad_{\xi_4}\xi_3},\xi_2}\xi_1+\\
        & ad_{\xi_4}ad_{ad_{\xi_3}\xi_2}\xi_1+ ad_{\xi_3} ad_{ad_{\xi_4}\xi_2}\xi_1 \biggm).
\end{split}
\end{equation}
}
{
Now we are in a position to provide the condition for the asymptotic equality in \eqref{eq:flow_equality} in light of the above series.
\begin{proposition}
\label{prop:asympequality}
    If 
    $ \tau= \frac{1}{\epsilon^{\hat{p}}}$ with $0<\hat{p}< 1$.
             Then,
\begin{equation}
\begin{split}
&\BigVert  \Bigl( \overrightarrow{exp} \int_0^\tau \epsilon \bm{f}_p dp - exp(\bm{V}_\tau^{(r)}) \Bigr) \phi \BigVert _{s,\mathbb{K}} \le \\
        & C_1 \Bigvert  \int_0^\tau \Vert \epsilon \bm{f}_p \Vert _{s+2r+2} dp \Bigvert ^{r+1} \Vert \phi \Vert _{s+r+2, O_{C_2}}=O(\epsilon^{r+1}\tau^{r+1}).
\end{split}
\end{equation}
    where the constants $C_1,C_2,O_{C2}$ depend on $s,r$, and $\mathbb{K}$.
\end{proposition}
}
\begin{proof}By \xref{prop:orderepsilontau}, $\int_0^\tau \Vert \epsilon \bm{f}_p \Vert_{r,\mathbb{K}'} dp \leq M\epsilon \tau\leq M\epsilon^{1-\hat{p}}$. Since $M$ is finite, and $\epsilon^{1-\hat{p}}\rightarrow 0$ as $\epsilon \rightarrow 0$, then there exists $\epsilon^{*}$ such that for all $\epsilon < \epsilon ^{*}$, $M\epsilon^{1-\hat{p}}\leq 1$. This
satisfies \cite[Proposition 4.1] {agrachev1978exponential}.
\end{proof}
\begin{remark}
{
   \xref{prop:asympequality} provides the error bound due to $r$-truncation in the flow expansion; this holds for both time-varying flow and autonomous flow since the right chronological exponential recovers the
regular exponential for autonomous vector fields.  \xref{prop:asympequality} provides the approximation condition between the trajectories generated via the time-varying flow (original system) and its autonomous
$r$-order approximation. This holds regardless of stability; a similar observation was made in \cite{vela2003general}.}
\end{remark}

{In a linear Floquet-like fashion, we move into finding a special case of an autonomous system that has its trajectory coinciding with the original time-periodic system, not at the unit time, but after a period $\tau={T}'$. 
 }    

\begin{theorem}[{Nonlinear Floquet Theorem}] \label{thm:nonlinearFloquet_prem}
   (\cite[Theorem 3.2]{sarychev2001lie},\cite[Theorem 5]{vela2003general})
    The flow $\stmboldsymbol{\bPhi}_\tau^{\epsilon \bm{f}}$ generated by the nonlinear time-periodic system in \eqref{eqn:dxdtauepsilon} can be represented as a composition $\stmboldsymbol{\bPhi}_\tau^{\epsilon \bm{f}} = e^{\stmboldsymbol{\bLambda} \tau} \circ \stmboldsymbol{\mathbbmbf{P}}_\tau$ of the flow $e^{\stmboldsymbol{\bLambda} \tau}$ of the autonomous vector field $\stmboldsymbol{\bLambda}$ and a $T'-$ periodic map $\stmboldsymbol{\mathbbmbf{P}}_\tau$ if the diffeomorphism $\stmboldsymbol{\bPhi}_{T'}^{\epsilon \bm{f}}$ admits a logarithm $\bm{V}_{T'}$, where
    {
\begin{equation}
\label{eq:Lambda_theorem}
    \stmboldsymbol{\bLambda} = \frac{1}{T'} \bm{V}_{T'} = \frac{1}{T'}\stmboldsymbol{\bPhi}_{T'}^{\epsilon \bm{f}}= \frac{1}{T'} {ln}\Bigl(\overrightarrow{exp}\Bigl(\int_0^{T'} \epsilon \bm{f}_\tau d\tau\Bigr)\Bigr).
\end{equation}
}
\end{theorem}
{\cite{sarychev2001lie} noted that by defining $\stmboldsymbol{\bLambda}$ as $\bm{V}_{T'}$ divided by the period $T'$ makes its flow $e^{\stmboldsymbol{\bLambda}\tau}$ generates averaged trajectories with  \xref{thm:nonlinearFloquet_prem}
being a nonlinear analog to the known linear Floquet theorem. The reader may refer to \cite[Table 1]{Maggia2020higherOrderAvg} for a nice comparison between linear and nonlinear Floquet theorems. As presented
in \cite{sarychev2001lie,vela2003general,Maggia2020higherOrderAvg}, the following system is the complete averaged system:}
\begin{equation}
\label{eqn:averaged_prem}
    \frac{d\bar{\bm{x}}}{d\tau} = \stmboldsymbol{\bLambda}(\bar{\bm{x}}) = \sum_{i=1}^r \stmboldsymbol{\bLambda}_i (\bar{\bm{x}}).
\end{equation}

{In \eqref{eqn:averaged_prem}, if $r= \infty$, then we have complete averaging. If $r=r^{*}$ is finite, then we have averaging of $r$-order.} In addition, the relation between the stability property of \eqref{eqn:averaged_prem} to that of the original nonlinear time-periodic system \eqref{eqn:dxdtauepsilon} is provided by the following theorem.
\begin{theorem}\label{thm:stability_prem} 
    ({\cite[Theorem 3.1]{sarychev2001lie}, \cite[Theorem 6]{vela2003general}}) If the monodromy map $\stmboldsymbol{\mathscrbold{M}} = \stmboldsymbol{\bPhi}_{T'}^{\epsilon \bm{f}}$ of the system \eqref{eqn:dxdtauepsilon} has a fixed point, then the flow $\stmboldsymbol{\bPhi}_{\tau}^{\epsilon \bm{f}}$, has a periodic orbit whose stability is determined by the stability of $\stmboldsymbol{\mathscrbold{M}}$.
\end{theorem}
\begin{remark}
{
The following can be concluded based on  \xref{thm:stability_prem} and statements available in references.
\begin{itemize}
  \item As clearly stated in \cite{vela2003general},  \xref{thm:stability_prem} means that if the complete (i.e.,~$r= \infty$) averaged autonomous system \eqref{eqn:averaged_prem} is asymptotically (exponentially)
stable for an equilibrium point, then the associated orbit of the original time-varying system is asymptotically (exponentially) stable. 
  \item If one concludes asymptotic stability for \eqref{eqn:averaged_prem} via a finite truncation $r=r^{*}$, then for  \xref{thm:stability_prem} to hold, one has to show or assume that the asymptotic stability
property for \eqref{eqn:averaged_prem} will not be distorted by any $r>r^{*}$. This is stated clearly in \cite[Section 8]{sarychev2001lie}.
  \item Given that $r=1$ in \eqref{eqn:averaged_prem} is equivalent to the first-order averaging in classical averaging, the results of \cite[Theorem 10.5]{khalil2002nonlinear} applies where exponential
stability of \eqref{eqn:averaged_prem} when $r=1$ implies the exponential stability of the original system; it is not proved, but we expect this to carry over to higher orders ($r>1$).  
\end{itemize}
}
\end{remark}
{As noted in \xref{TeXFolio:rmk4} and \xref{TeXFolio:rmk5}, the approximation between the trajectories of \eqref{eqn:dxdtauepsilon} and the averaged system trajectories \eqref{eqn:averaged_prem} for finite $r=r^{*}$ is established. {Now, we prove that the bound in  \xref{prop:asympequality}, i.e.,~$O(\epsilon^{r+1}\tau^{r+1})$ can be made arbitrary small for any time horizon that can be stretched as needed.}
{
\begin{proposition}
\label{prop:vanishingBounds}
              Let the time horizon $\tau=1/\epsilon^{\tilde{p}}$ where $0<\tilde{p}<1$. Then, the bounds in  \xref{prop:orderepsilontau,prop:asympequality} vanish as $\epsilon\rightarrow 0$ (or equivalently $\omega \rightarrow \infty$). 
\end{proposition}
\begin{proof}It is obvious that as $\epsilon \rightarrow 0$, $\tau=1/\epsilon^{\tilde{p}} \rightarrow \infty$. 
  From  \xref{prop:orderepsilontau}, $\int_0^\tau \Vert \epsilon \bm{f}_p \Vert_{r,K'} dp \leq M\epsilon \tau$ for some finite $M>0$. Observe that  
 $\epsilon \tau= \epsilon/\epsilon^{\tilde{p}}= \epsilon^{1-\tilde{p}}$; in terms of $\omega$, we substitute $\epsilon= \omega^{p^{*}-1}$ as it is defined earlier. This results in $\epsilon \tau= \epsilon^{1-\tilde{p}}= \omega^{(p^{*}-1)(1-\tilde{p})}$. Clearly with $0<p^{*}<1$ and $0<\tilde{p}<1$, it follows that $\epsilon\tau \rightarrow 0$ as $\epsilon \rightarrow 0$ (or equivalently as $\omega \rightarrow \infty$). Since $r\geq1$, the bound in  \xref{prop:asympequality} vanishes as well (i.e.,~$\epsilon^{r+1}\tau^{r+1} \rightarrow 0$) as a result of $\epsilon\tau \rightarrow 0$. 
\end{proof}
}
}

{
Moreover, one can directly relate the trajectories approximation established above, asymptotic stability of the averaged system of order $r=r^{*}$, and the concept of practical stability \cite{moreau2000practical}
straightforwardly. We provide the following theorem.}
\begin{theorem}
\label{thm:rstar}
{
         If the system \eqref{eqn:averaged_prem} for some $r=r^{*}$ has an equilibrium point $\bar{\bm{x}}^{*}\in \mathbb{D}_0$ that is asymptotically locally uniformly stable, then the system \eqref{eqn:dxdtauepsilon} is practically uniformly asymptotically stable for $\mathbb{D}_0$.    
    }
\end{theorem}
{
\begin{proof}
It follows immediately from  \xref{prop:vanishingBounds} that for any fixed time horizon $\tau_f>0$, the error bound between the solution of \eqref{eqn:dxdtauepsilon} and solution of \eqref{eqn:averaged_prem} in the time domain $\tau \in [0,\tau_f]$ can be made arbitrarily small. That is, for any $\tau_f>0$, $\epsilon^{r+1}\tau_f^{r+1} \rightarrow 0$ as $\epsilon \rightarrow 0$, or equivalently as $\omega \rightarrow \infty$.
\end{proof}
}

{
Finally, for computing $r$-order averaging (also referred to in literature sometimes as logarithms) for the averaged system in \eqref{eqn:averaged_prem}, one should compute $r$ terms (partial sum) of the
series \eqref{eqn:voltera}. Simplifications using properties of Lie brackets can be used so that one can get $\mathscr{G}_r$ expressions (e.g.,~\eqref{eqn:commutators}). After that, one divides by $T'$ (per
\eqref{eq:Lambda_theorem}). For more details, the reader can refer \cite[Section 2.2.2]{vela2003averagingThesis}. First four logarithms are provided below \cite{Maggia2020higherOrderAvg}:
      }
\begin{align}
&\stmboldsymbol{\bLambda}_1 = \frac{\epsilon}{T'}\int_0^{T'} \bm{f}_\tau d\tau, \label{eqn:lambda1_main} \\
&\stmboldsymbol{\bLambda} _2  = \frac{\epsilon^2}{2!T'}\int_0^{T'} \left [\int_0^\tau \bm{f}_p dp, \bm{f}_\tau\right ] d\tau, \label{eqn:lambda2_main}\\
&\begin{array}{ll}
\stmboldsymbol{\bLambda}_3 &= \frac{\epsilon^3}{3!}\Biggm( -\frac{3}{2}[\stmboldsymbol{\bLambda}_1, \stmboldsymbol{\bLambda}_2] +\\
        &\frac{2}{T'} \int_0^{T'} \left [\int_0^\tau \bm{f}_p dp, \left[\int_0^\tau \bm{f}_p dp, \bm{f}_\tau \right ] \right] d\tau \Biggm)
\end{array}\label{eqn:lambda3_main}\\
&\begin{array}{ll}
\stmboldsymbol{\bLambda}_4 &= \frac{2 \epsilon^4}{4!T'} \int_0^{T'}  \biggm( \int_0^\tau \left[ \int_0^p \left[\int_0^q \bm{f}_r dr, \bm{f}_q \right]dq, [\bm{f}_p,\bm{f}_\tau] \right] dp\\
&+\xs  \left[ \int_0^\tau \left[ \int_0^p \left[ \int_0^q \bm{f}_r dr , \bm{f}_q\right] dq, \bm{f}_p\right] dp,\bm{f}_\tau \right]\\
&+\xs  \int_0^\tau \left[ \int_0^p \bm{f}_q dq, \left[\left[ \int_0^p \bm{f}_q dq, \bm{f}_p\right],\bm{f}_\tau\right]\right]dp \biggm)d\tau.
\end{array}\label{eqn:lambda4_main}
\end{align}

\section{Higher-order {\?{Lie}} bracket approximation and application to extremum seeking}

\begin{theorem}
\label{thm:mainthm_main}
Let Assumptions A1-A2, the results of  \xref{lemma:wellposed_main} and \xref{TeXFolio:pps4} hold, then for $r=1,\ldots,4$, the ($r$)-order averaging of the control system in \eqref{eqn:controlAffineIntro} corresponds directly to the $(r-1)$-order LBS provided in \eqref{eqn:LBS_main}:
\begin{equation}
\label{eqn:LBS_main}
    \dot{{\bm{z}}} = \bm{b}_{\bm{0}}({\bm{z}})+ \sum_{i = 1}^r \bm{L}_i({\bm{z}}),
\end{equation}
with
\begin{align}
\bm{L}_1 &= 0, \label{eqn:l1_main}\\
        \bm{L}_2 &= \sum_{j_1 = 1}^m \sum_{j_2=j_1+1}^m \nu_{j_1 j_2}[\bm{b}_{j_1}, \bm{b}_{j_2}],\label{eqn:l2_main}\\
        \bm{L}_3 &= \sum_{j_1 = 1}^m \sum_{j_2=j_1+1}^m \sum_{j_3=1}^m  \nu_{j_1 j_2 j_3} [\bm{b}_{j_3},[\bm{b}_{j_1},\bm{b}_{j_2}]],\label{eqn:l3_main}\\
\label{eqn:l4_main}
\nf&\begin{array}{ll}
\bm{L}_4 &= \displaystyle \sum_{j_1 = 1}^m \sum_{j_2=j_1+1}^m \sum_{j_3=1}^m \sum_{j_4=j_3+1}^m \beta_{1_{j_1 j_2 j_3 j_4}} \Bigl[[\bm{b}_{j_1},\bm{b}_{j_2}],[\bm{b}_{j_3},\bm{b}_{j_4}]\Bigr]\\
&+ \displaystyle  \sum_{j_1 = 1}^m \sum_{j_2=j_1+1}^m \sum_{j_3=1}^m \sum_{j_4=1}^m \beta_{2_{j_1 j_2 j_3 j_4}} [[[\bm{b}_{j_1},\bm{b}_{j_2}],\bm{b}_{j_3}],\bm{b}_{j_4}],
\end{array}
\end{align}
    where
\begin{equation}
\begin{split}
\nu_{j_1 j_2} &= \frac{\omega^{p_{j_1}+p_{j_2}-1}}{2T'} \int_0^{T'}  \Bigm( u_{j_2}(k_{j_2}\tau)\int_0^\tau u_{j_1}(k_{j_1} p)dp \\
&-\xs  u_{j_1}(k_{j_1}\tau)\int_0^\tau u_{j_2}(k_{j_2} p)dp \Bigm)d\tau,
\end{split}\label{eqn:nuij_main}
\end{equation}
\xpar\cvskip[-9pt]
\begin{equation}
\begin{split}
\nu_{j_1 j_2 j_3} &= \frac{\omega^{p_{j_1}+p_{j_2}+p_{j_3}-2}}{3T'}\int_0^{T'} \int_0^\tau u_{j_3}(k_{j_3} s)ds \Bigm(u_{j_2}(k_{j_2} \tau )\\
        &\int_0^\tau u_{j_1}(k_{j_1} p) dp 
         -u_{j_1}(k_{j_1}\tau)\int_0^\tau u_{j_2}(k_{j_2} p)dp \Bigm) d \tau,
\end{split}\label{eqn:nuijk_main}
\end{equation}
    and
\begin{align*}
\beta_{1_{j_1 j_2 j_3 j_4}} &= \omega^{p_{j_1}+p_{j_2}+p_{j_3}+p_{j_4}-3} \frac{1}{12T'}\int_0^{T'}  \alpha_5(\tau)_{j_1 j_2 j_3 j_4}d\tau\\
        \beta_{2_{j_1 j_2 j_3 j_4}} &= \omega^{p_{j_1}+p_{j_2}+p_{j_3}+p_{j_4}-3} \frac{1}{12T'}\int_0^{T'} \bigl(\alpha_8(\tau)_{j_1 j_2 j_3 j_4} \\
&-\xs  \alpha_{10}(\tau)_{j_1 j_2 j_3 j_4}\bigr)d\tau,
\end{align*} 
    with
\begin{align}
\nf&\begin{array}{ll}
\alpha_5(\tau)_{j_1 j_2 j_3 j_4} &= \int_0^\tau \biggm( \Bigl(u_{j_4}(k_{j_4}\tau) u_{j_3}(k_{j_3}p) - u_{j_3}(k_{j_3}\tau) u_{j_4}(k_{j_4}p) \Bigr) \mrowclose\\
         & \times\Bigm( \int_0^p \bigl (u_{j_2}(k_{j_2}q)
          \mrowopen\int_0^q u_{j_1}(k_{j_1} r)dr -\mrowclose\\
          & u_{j_1}(k_{j_1}q)\int_0^q u_{j_2}(k_{j_2} r)dr \bigr)dq \Bigm) \biggm)dp,
\end{array}
\\
\nf&\begin{array}{ll}
\alpha_8(\tau)_{j_1 j_2 j_3 j_4} &= u_{j_4}(k_{j_4}\tau) \int_0^\tau \biggm( \int_0^p \Bigm( u_{j_2}(k_{j_2}q)\int_0^q u_{j_1}(k_{j_1} r)dr- \\
        & u_{j_1}(k_{j_1}q)\int_0^q u_{j_2}(k_{j_2} r)dr \Bigm) dq \, u_{j_3}(k_{j_3}p)  \biggm) dp ,
\end{array}
\\
\nf&\begin{array}{ll}
\alpha_{10}(\tau)_{j_1 j_2 j_3 j_4} &= \int_0^\tau \biggm(  \Bigl( u_{j_2}(k_{j_2}p)\int_0^p u_{j_1}(k_{j_1} q)dq-u_{j_1}(k_{j_1}p) \\
& \int_0^p u_{j_2}(k_{j_2} q)dq \Bigr)
        u_{j_3}(k_{j_3}\tau)
         \int_0^p u_{j_4}(k_{j_4} q)dq \biggm) dp;
\end{array}
\end{align}
    where $p,q,r,s,\tau$ are the variables of integration, and $T' = (2\pi)\times LCM(k_1^{-1}, k_2^{-1}, \ldots ,k_m^{-1})$.
\end{theorem}
\textit{Proof of \xref{thm:mainthm_main} is provided in  Appendix A.}

\begin{corollary}
\label{cor:bounded_main}
    If $\bm{z}(t) \in \mathbb{D}_0$, for any $t_f>0$, $\forall t \in [0,t_f]$, then there exists $\epsilon^{*}>0$ such that for all $\epsilon\in(0,\epsilon^{*})$, there exists $d(\epsilon)>0$ such that {$\svert \bm{x}(t)-\bm{z}(t)\svert  \leq d(\epsilon)$ for $t \in [0,t_f]$}. Moreover, as $\epsilon \rightarrow 0$, $d(\epsilon) \rightarrow 0$ for any $r$.   
\end{corollary}
\begin{remark}
\label{rem:bounded_main}
     \xref{cor:bounded_main} follows immediately from  \xref{thm:mainthm_main} where it is shown that the trajectories of the higher-order LBS denoted by $\bm{z}$ are in fact trajectories of a higher-order averaged system; hence,  \xref{prop:vanishingBounds} applies. 
 In addition,  \xref{cor:bounded_main} generalizes what needed to be proved in \cite[Theorem 1]{DURR2013} and \cite[Lemma 1]{labar2019newton}, where the distance between the original system and LBS is shown to be
bounded.   
         
\end{remark}

\begin{corollary}
\label{cor:stability_main}
{Suppose for finite $r=r^{*}$ we have $\bm{z}^{*}\in \mathbb{D}_0$ an equilibrium point of the system \eqref{eqn:LBS_main} that is locally uniformly asymptotically stable. Then \eqref{eqn:controlAffineIntro} is locally practically uniformly asymptotically stable for $\mathbb{D}_0$.}  
\end{corollary}
    
\begin{remark}
\label{rem:stability_main}
     \xref{cor:stability_main} is a {straightforward result from  \xref{thm:mainthm_main} where  \xref{cor:bounded_main} and  \xref{thm:rstar} apply.}
    
\end{remark}

\subsection{Complete averaging by design, asymptotic stability results, and re-examination of literature works}

{From  \xref{prop:orderepsilon} and for finite $T'$ and $r$, \eqref{eqn:averaged_prem} can be written as (similar observation is made in \cite{vela2003general}):
\begin{equation}
\label{eqn:completer_averaging_inTau_epsilon_rPLUS1}
    \frac{d\bar{\bm{x}}}{d\tau } = \sum _{i=1}^r \Lambda_i+O(\epsilon ^{r+1}).
\end{equation}
In $t-$scale, another representation of \eqref{eqn:LBS_main} by \eqref{eqn:completer_averaging_inTau_epsilon_rPLUS1} is:
\begin{equation}
\dot{\bm{z}} = \sum _{i=1}^r \bm{L}_i(\bm{z},\omega)+ O(\epsilon^{r+1} \omega).
\end{equation}
}
\begin{proposition}
\label{prop:complete_avg}
{
    If $\lim_{\omega \rightarrow \infty} \sum_{i=1}^r \bm{L}_i(\bm{z},\omega) = \bm{G}(\bm{z})<\infty, \; \forall \bm{z} \in \mathbb{D}_0$ and $0<max(p_i) \le p^{*} < r/(r+1)$, then $\dot{\bm{z}}=\bm{G}(\bm{z})$ represents a complete averaging asymptote.}
\end{proposition}
\begin{proof}
We want to find the condition for $\epsilon^{r+1} \omega \rightarrow 0$ as $\omega \rightarrow \infty$. Note that $\epsilon^{r+1} \omega = \omega^{(p^{*}-1)(r+1)+1}$, so we need $(p^{*}-1)(r+1)+1<0$. This is satisfied if $p^{*}r+p^{*}-r<0$, i.e.~$p^{*}<r/(r+1)$.
\end{proof}
\begin{remark}
\label{rem:complete_avg} { \xref{prop:complete_avg} provides
a condition that can be met by design for which one can make an achievable choice of $max\{p_i\}$, i.e.,~all $p_i<r/(r+1)$ in \eqref{eqn:controlAffineIntro} so that the \textit{ultimate} behavior of
\eqref{eqn:controlAffineIntro} is qualitatively captured by $\dot{\bm{z}}=\bm{G}(\bm{z})$ as $\omega \rightarrow \infty$. This also means that  \xref{thm:stability_prem} applies. That is, asymptotic (exponential) stability of $\dot{\bm{z}}=\bm{G}(\bm{z})$ implies
asymptotic (exponential) stability of \eqref{eqn:controlAffineIntro} as $\omega \rightarrow \infty$ (in the sense of \cite[Theorem 10.5]{khalil2002nonlinear}).}
    
\end{remark}
{The condition provided in  \xref{prop:complete_avg} is also explainable from our generalized LBS in \eqref{eqn:LBS_main}. For instance, if $r=2$ (1st-order LBS), then $p_1$ and $p_2$ will be both
strictly less than $r/(r+1)=2/3$. Hence, going higher in order, the power of $\omega$ will be strictly negative since any {$p_{j_1}+p_{j_2}+p_{j_3}<2/3+2/3+2/3=2$}, hence $\nu_{j_1 j_2 j_3}$ will vanish as $\omega \rightarrow \infty$; it remains to show that
$\lim_{\omega \rightarrow \infty}  L_2(\bm{z},w)=G(\bm{z})$ bounded and defined. Similarly, if $r=3$, 2nd-order LBS, all $p_i$ less than $3/4$ means that $p_{j_1}+p_{j_2}+p_{j_3}+p_{j_4}-3<0$ and $\beta_{1_{j_1 j_2 j_3 j_4}}$, $\beta_{2_{j_1 j_2 j_3 j_4}}$ vanish as $\omega \rightarrow \infty$; it remains then that one
shows $\lim_{\omega \rightarrow \infty}(  \bm{L}_2(\bm{z},w)+\bm{L}_3(\bm{z},w))=\bm{G}(\bm{z})$ bounded and defined. Now, let us re-examine the works of \cite{DURR2013,BoundedUpdateKrstic,VectorFieldGRUSHKOVSKAYA2018}.}
\begin{corollary}
\label{cor:complete_avg}
    {If $p_i=1/2$ for all $i$ in \eqref{eqn:controlAffineIntro}, then $\dot{\bm{z}}= \bm{L}_2$ in \eqref{eqn:LBS_main} represents complete averaging asymptote.}
\end{corollary}
\begin{remark}
{ \xref{cor:complete_avg} can be easily verified by examining $\nu_{j_1 j_2 j_3}$. Clearly any $p_{j_1}+p_{j_2}-1=0$, hence $\lim_{\omega \rightarrow \infty}  \bm{L}_2$ is bounded, satisfying  \xref{prop:complete_avg}. This also means that  \xref{rem:complete_avg}
applies for asymptotic stability. To conclude, in systems like those in \cite{DURR2013,BoundedUpdateKrstic,VectorFieldGRUSHKOVSKAYA2018}, first-order LBS approximation suffices for capturing the ultimate qualitative
behavior of the system for stability, and convergence rate, among other qualitative behaviors.}  
\end{remark}
{Now we re-examine the work of \cite{labar2019newton}. Going to higher orders such as $r=3$ is more challenging even if one chooses for example $p_i=2/3$ for all $i$ to satisfy $max\{p_i\}<3/4$ per
\xref{prop:complete_avg} and also to maintain the power of $\omega$ at 0 in $\nu_{j_1 j_2 j_3}$. That is, $p_i=2/3$ for all $i$ implies that $p_{j_1}+p_{j_2}>1$, making the power of $\omega$ greater than 0 in
$\nu_{j_1 j_2}$; this makes the limit condition in  \xref{prop:complete_avg} \textit{not applicable} as $\nu_{j_1 j_2}$ is unbounded when taking limit $\omega \rightarrow \infty$. This reality imposes a challenge to acquire a complete
averaging asymptote as in  \xref{prop:complete_avg}. One has to involve the Lie brackets themselves in $\bm{L}_2$ in \eqref{eqn:LBS_main} and/or the iterated integrals (taking advantage of their different frequency
scales, $k_i$); hence \cite[Assumption 2]{labar2019newton}. That assumption imposed that the Lie brackets and/or iterated integrals in all undesired terms in $\bm{L}_2$ and $\bm{L}_3$ have to vanish whenever
the power of $\omega$ is positive; we only consider these terms with power of $\omega$ at 0 as desired terms in the design. The mentioned assumption is only a special case. One can surely keep powers of
$\omega$ in undesired terms positive and still maintain a LBS approximation for fixed $\omega$ per averaging theory. 
 Note that our results \eqref{eqn:LBS_main} when we truncate at $r=3$ recover the second-order LBS by \cite{labar2019newton} but without the imposed $lim \omega \rightarrow \infty$ (see \eqref{eqn:Lie_Prem}) as this is a special
case of our results} by doing the following with \eqref{eqn:nuijk_main}: (i) change the variable $\tau$ to $t$ using the relation $\tau = \omega t$; (ii) impose that $u_{j_1}(k_{j_1}\omega\tau)\int_0^\tau u_{j_2}(k_{j_2} \omega p) = -u_{j_2}(k_{j_2}\omega \tau)\int_0^\tau u_{j_1}(k_{j_1} \omega p)$ (e.g.,~specific kinds of
sinusoidal functions). {The second-order LBS as it was used in ESC design by \cite{labar2019newton} is a complete averaging asymptote by the following more relaxed condition derived from  \xref{prop:complete_avg} which is
important in relating not only $p_i$ with $r$, but also taking in consideration the number of control inputs $m$.}       

\begin{corollary}
\label{cor:completeAvg_main}
Let $m<r$, $\lim_{\omega \rightarrow \infty}\bm{L}_i$ bounded and defined for $i=1,\ldots,m$, and $p_{j_1}+p_{j_2}+\cdots+p_{j_m}=m-1$, then with $\omega \rightarrow \infty$, $\bm{L}_{m+1}= \bm{L}_{m+2}= \cdots= \bm{L}_r = 0$ and
\begin{equation}
\dot{{\bm{z}}} = \bm{b}_{\bm{0}}({\bm{z}})+ \sum_{i = 1}^m \bm{L}_i({\bm{z}}),
\end{equation}
    provides a complete averaging asymptote.
\end{corollary}

\begin{remark}
\label{rem:complete_avg_main}
     \xref{cor:completeAvg_main} means that if for $m$ control inputs, if one chooses $p_{j_1}+p_{j_2}+\cdots+p_{j_m}=m-1$, then as $\omega \rightarrow \infty$, higher-order terms $\bm{L}_{m+i}$ will {vanish due to negative powers of $\omega$}, leaving the remaining terms ($\bm{L}_{1}, \ldots ,\bm{L}_{m}$) to capture the completed averaged properties of the original control-affine system. For example, if $m=2$, and $p_1=p_2=0.5$, i.e.,~$p_1+p_2=1$, then one can see that the power of $\omega$ will be negative for $\bm{L}_3$ {because adding $p_1$ or $p_2$ (both of them are strictly less than 1) to $p_1+p_2$ will be less than 2}(see \eqref{eqn:l3_main}) and similar observation apply to $\bm{L}_4$ (see \eqref{eqn:l4_main}).
\end{remark}

\subsection{Application to extremum seeking and simulations}

 With the understanding of  \xref{thm:mainthm_main} which establishes the connection between higher-order averaging and LBS approximations for control-affine ESC systems, and along with the strong insights derived from  \xref[env={Corollaries}]{cor:bounded_main} to \nxref{cor:completeAvg_main} and {\xref{TeXFolio:pps5}}, one can have a more comprehensive understanding of control-affine ESC systems and how to design them for better performance {and based on what order of LBS}. This perhaps will need many studies to do after the results this paper provides, but here we will highlight a few of the observations we were able to make and verify them via Examples.
\begin{enumerate}[label=$O$\arabic*.]
    \item {For a given ESC, one can check the condition of \xref{TeXFolio:pps5} (\xref{TeXFolio:rmk8}) or \xref{TeXFolio:cry4} (\xref{TeXFolio:rmk10}) to determine the r-order of LBS that is complete averaging asymptote which is the ultimate representative of the ESC qualitative behavior (e.g.,~convergence rate and stability). See Examples \xref[range]{TeXFolio:exm1,TeXFolio:exm2}.}
    \item {Thanks to \eqref{eqn:LBS_main}, one can improve ESC performance by design (especially convergence rate) to be influenced by   higher-order derivatives, even beyond Hessian, which is unprecedented in ESC literature to the best of our knowledge, while maintaining a bounded average/LBS as $\omega \rightarrow \infty$. In Examples \xref[range]{TeXFolio:exm3,TeXFolio:exm4} we show how one by careful choice of $p_i$, $k_i$ with tracking zero and non-zero Lie brackets, can achieve an ESC with a bounded average as $\omega \rightarrow \infty$ that estimates third-order derivative information for the first time.}
         \item  Due to higher-order averaging theory, truncating before the desired order may or may not capture stability, but even if it captures stability, convergence behavior will best be represented by the complete averaging asymptote. Refer to Examples \xref[range]{TeXFolio:exm2,TeXFolio:exm3}.
  \item {ESC designer can totally make choices that do not lead to complete averaging asymptote. LBS approximations in \eqref{eqn:LBS_main} still hold and provide the results of \xref[range]{TeXFolio:cry1,TeXFolio:cry2}; no need to impose
$lim\:  \omega \rightarrow \infty$ as in \cite{labar2019newton} (also provided in \eqref{eqn:Lie_Prem})}.
\end{enumerate}

\begin{figure}
\includegraphics{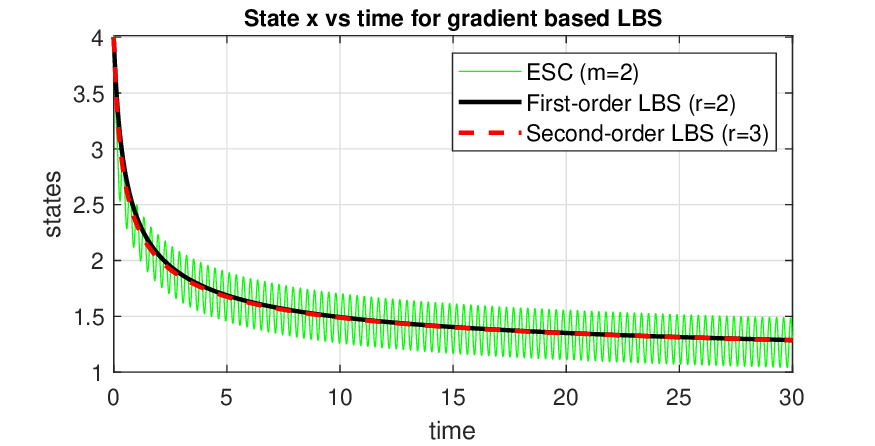}
\caption{\label{fig:2contorls2lbs1_sim}Comparison of ESC with 2 control inputs and LBS approximations truncated at $r=2$ and $r = 3$.}
\end{figure}

\textbf{Example 1.} This example verifies observation $O1$ that {if \xref{TeXFolio:pps5} is satisfied (like the works of \cite{DURR2013,BoundedUpdateKrstic,VectorFieldGRUSHKOVSKAYA2018}
as in \xref{TeXFolio:cry3}), then we have a complete averaging asymptote} that captures the ultimate qualitative behavior of the ESC (e.g.,~stability and convergence rate), and higher orders will not be needed. Let us consider:
\begin{equation}
\label{eqn:example1_sim1}
\dot{x}= \sum_{i=1}^{m} \omega^{p_i} b_i(x) u_i(k_i\omega t),
    \dot{v} = h(J-v);
\end{equation}
with $m = 2, J=-H(x-1)^4$ as the cost function and $v$ as a variable introduced by the presence of a high pass filter of frequency $h$. The expressions and numerical values used in our simulation are $b_1 = J-v$, $b_2 = a$, $u_1 = \sin(k_1 \omega t), u_2 = \cos(k_2 \omega t)$, $p_1=p_2= 0.5$, $k_1=k_2=1$, $\omega = 20$, $h=5$, $a=1$ and $H=1/10$. The initial condition is $[x_0,v_0]=[4,0]$. We computed first-order LBS (r=2) and second-order LBS (r=3) using the formulas of $\bm{L}_1$ to $\bm{L}_3$ in \eqref{eqn:l1_main}--\eqref{eqn:l3_main}, and plotted both of them with the ESC system.
 The simulation is shown in  \xref{fig:2contorls2lbs1_sim}. One can clearly see how well the LBSs capture the average behavior of the ESC system. It is also clear that the first-order LBS (r=2) was enough (no need for LBS with r=3) to capture the qualitative behavior of the system (e.g.,~stability and convergence). Furthermore, note that the number of control inputs is $m=2$ with $p_1+p_2 = 1$. So,  \xref{cor:completeAvg_main} {is applicable as well}.

\begin{figure}
\includegraphics{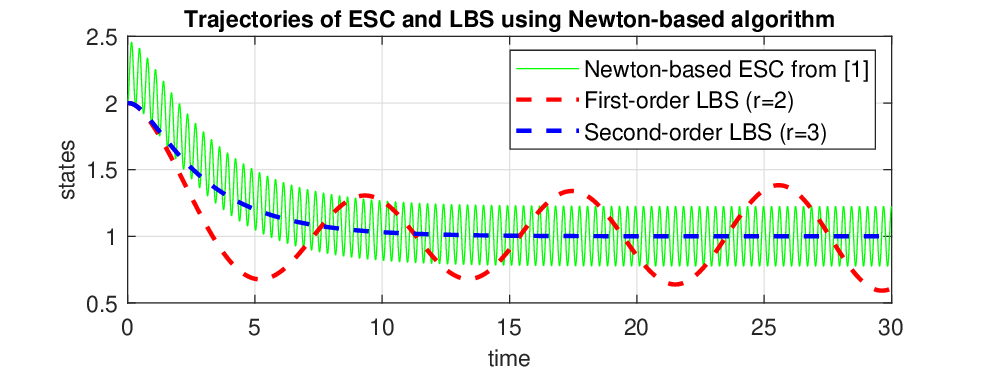}
\caption{\label{fig:NewtonLbs_sim}Trajectories of ESC and LBS using Newton-based algorithm with LBS truncated at $r=3$ (dashed-blue) and $r=2$ (dashed-red).}
\end{figure}

\textbf{Example 2.} In this example we demonstrate {two things. First, the importance and need in some cases for going to higher-order LBSs to capture stability and achieve complete averaging asymptote as in observation $O1$. Second, if the design needs higher order derivative information such as Hessian, then as in $O2$ we may need at least $m=3$ so second-order Lie brackets extract Hessian from vector fields involving the objective function.}
We show that taking $r=2$ in \eqref{eqn:LBS_main} to approximate the Newton-based ESC in \cite{labar2019newton} will not be sufficient to capture the averaged dynamics of the ESC system and one needs to go higher,
to $r=3$, to capture the qualitative behavior of the ESC {via a complete averaging asymptote}. We take the system in \cite[Section 5]{labar2019newton}. The ESC system is in the form of
\eqref{eqn:controlAffineIntro} with states $\bm{x} = [x_1,d,y,z]^T$, $\bm{b}_0 = [\rho d, -\omega_d(y+zd), -\omega_y y, \omega_z z]^T$, $\bm{b}_1 = [1,0,0,0]$, $\bm{b}_2 = [0,0,a_2J,0]$, $\bm{b}_3 = [0,0,0,a_3J]$ with $J = H(x_1-1)^2$, $a_2 = -2k\omega_y$, $a_3 = 8k^2 \omega_z$. The values of parameters are $H = 2, k=1, \rho = 0.3, \omega_d = 0.5, \omega_y = 20, \omega_z = 0.5$, $p_1 = 0.51, p_2 = 0.49, p_3 = 0.98$. The
control inputs are $u_1 = \sin(k\omega t)$, $u_2 = \cos(k\omega t)$ and $u_3 = \cos(2k\omega t)$ with $\omega = 20$. Now, the higher-order LBS can be calculated 
as $\bm{z} = \bm{b}_0+[0,0,\omega_y \nabla J(z),\omega_z \nabla^2J (z)]^T$. The simulation results are shown in  \xref{fig:NewtonLbs_sim}.
One can clearly see that the second-order LBS captures the qualitative behavior (e.g.,~stability and convergence) of the ESC system. However, if one truncates the LBS at the first-order ($r=2$), then the Hessian information is lost and the system actually becomes unstable, i.e.,~(it does not capture any qualitative behavior, e.g.,~stability or convergence).

\begin{figure}
\includegraphics{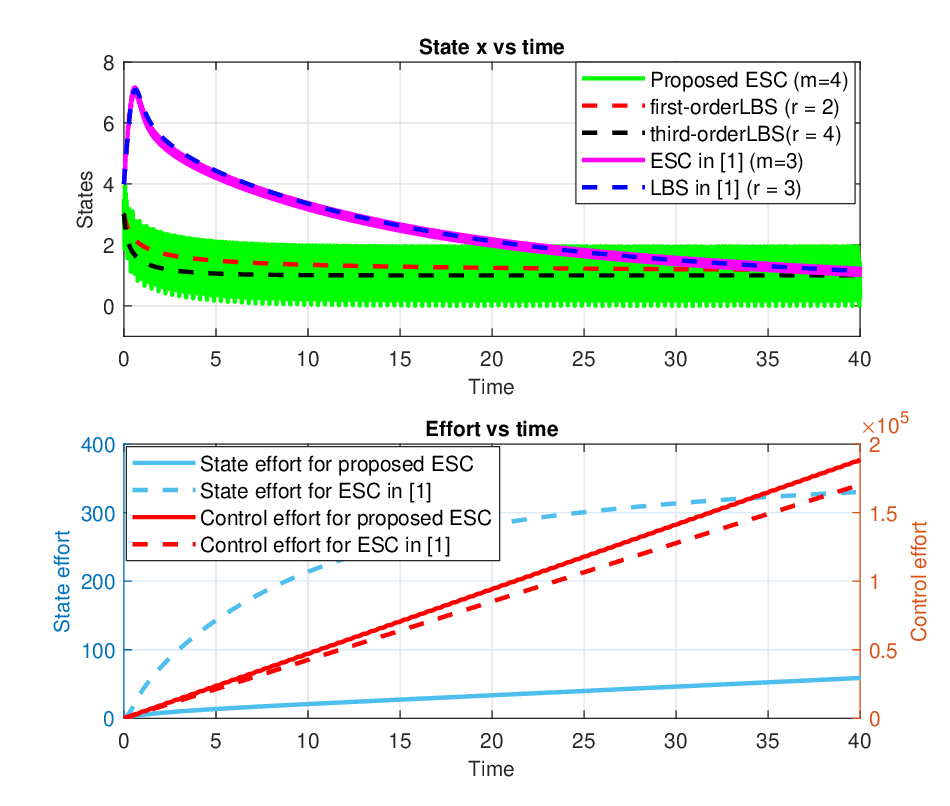}
\caption{\label{fig:4control3lbs_sim}Comparison of performance between our simple, proposed ESC with $m=4$ based on third-order LBS vs. Newton-based ESC with $m=3$ from \cite{labar2019newton}.}
\end{figure}

\textbf{Example 3.}
{The objective here as stated in $O2$ is to design an ESC with third-order derivative influence that possesses a high convergence rate. Let us consider system \eqref{eqn:example1_sim1} but with $m=4$ such
that $b_1 = J-v$, $b_2 = b_3 = b_4=a$, $u_1 = \sin(k_1 \omega t), u_2 = \cos(k_2 \omega t), u_3 = \sin(k_3 \omega t), u_4 = \cos(k_4 \omega t)$. We take a fourth-degree polynomial cost function $J=-H(x-1)^4$ to test convergence better. Our proposed ESC has $p_1=p_2=p_3 = 0.5, p_4 = 0.99$, $k_1=k_2=1, k_3 = 1/3, k_4 = 3/2$. For simulation we take $\omega = 100$,
$h=5$, $a=1$ and $H=1/5$. The initial condition is $[x_0,v_0]=[3,0]$. Here we briefly explain our design rationale, but the reader is directed to our supplementary file and code \cite{github} for all
details. Now, for $L_2$ corresponding to first-order Lie brackets (see \eqref{eqn:l2_main}), we have $[b_2,b_3]=[b_2,b_4]=[b_3,b_4]=0$. So, $[b_1,b_2],[b_1,b_3],[b_1,b_4]$ are the relevant non-zero Lie brackets (all of which depend on $[b_1,b_2]$); and
only $\nu_{12},\nu_{13},\nu_{14}$ are relevant, however, only $\nu_{12}$ has non-zero value given by $\nu_{12}=0.5$.
Now, for $L_3$ corresponding to the second-order Lie brackets (see \eqref{eqn:l3_main}), the Lie brackets involving $[b_2,b_3]$, $[b_2,b_4], \text{or } [b_3,b_4]$ will be zero and the relevant Lie brackets are given by $[b_1,[b_1,b_2]]$ and $[b_2,[b_1,b_2]]$. There are 12 relevant $\nu_{j_1j_2j_3}$ (corresponding to non-zero Lie brackets) given by $
        {\nu_{121}}, \nu_{122}, {\nu_{123}}, \bm{\nu}_{\bm{124}}, 
        {\nu_{131}}, \nu_{132}, {\nu_{133}}, \bm{\nu}_{\bm{134}} ,
        \bm{\nu}_{\bm{141}}, \bm{\nu}_{\bm{142}},\\ \bm{\nu}_{\bm{143}}, \bm{\nu}_{\bm{144}},  
$
where bold $\nu_{j_1j_2j_3}$ corresponds to those with $p_{j_1}+p_{j_2}+p_{j_3} \approx 2$ except for $\bm{\nu}_{\bm{144}}$, for which $p_{j_1}+p_{j_2}+p_{j_3} > 2$, threatening the boundedness of $L_3$ as $\omega \rightarrow \infty$. However, thanks to the choice of control inputs, the corresponding iterated integral (see \eqref{eqn:nuijk_main}) neutralized the effect of $\bm{\nu}_{\bm{144}}$. Furthermore, it is found that the iterated integrals are not all zeros, thus, making $L_3$ significantly important in the higher-order LBS calculation.
Now, for $L_4$ corresponding to the third-order Lie brackets (see \eqref{eqn:l4_main}), one will find that the Lie brackets of the form $[[b_{j_1},b_{j_2}],[b_{j_3},b_{j_4}]]$ are zero since $b_2 = b_3 = b_4$, and the non-zero relevant Lie brackets are $[[[b_1,b_2],b_1],b_1],[[[b_1,b_2],b_1],b_2],[[[b_1,b_2],b_2],b_1]$ and $[[[b_1,b_2],b_2],b_2]$. Since the Lie brackets multiplied with $\beta_{1_{j_1j_2j_3j_4}}$ (see \eqref{eqn:l4_main}) are all zeros, $\beta_{1_{j_1j_2j_3j_4}}$ do not need to be computed and there are only 48 relevant $\beta_{2_{j_1j_2j_3j_4}}$ that need to be computed. All of the $\beta_{2_{j_1j_2j_3j_4}}$ have not been provided due to space constraints but it is found that $\beta_{2_{j_1j_2j_3j_4}}$ corresponding to $p_{j_1}+p_{j_2}+p_{j_3}+p_{j_4} \approx 3$ do not all have zero iterated integrals (with some, actually, having significant values). Thus, $L_4$ (capturing third-order derivative information) has considerable importance in the higher-order LBS. Note that $\beta_{2_{1444}} = 0$, hence $L_3$ is bounded with respect to $\omega$, so is the entire third-order LBS for this ESC system.  
}
    We compare the performance of the proposed ESC with the Newton-based ESC and its LBS as in Example 2 but with the objective function and relevant parameters as the proposed system. The initial condition for this system is chosen as $\bm{x}_0 = [4,26.6,0,0]$ such that two systems will have the same initial convergence rate for performance comparison.
 \xref{fig:4control3lbs_sim} {(top)} shows the simulation results comparing both cases.

It is remarkable and immediately noticeable that the rate of convergence of the simple, proposed ESC is very high (converge within 5--10~s) and still qualitatively captured by the third-order LBS in a much better manner
than if one truncates at first-order LBS; however, the first-order LBS is still sufficient for providing the stability property similar to \cite{DURR2013}. On the other hand, the Newton-based ESC system (the second
system) is much slower (converges within 40~s). The reason for the faster convergence for the first system is 
 the extracted third-order derivative information influencing the complete averaging asymptote (i.e.,~the third-order LBS). 
  {For further analysis, we also present in  \xref{fig:4control3lbs_sim} {(bottom)} the effort exerted by the control inputs and the state effort of the system. To find the control effort we use the formula $\int_0^t \sum_{i=1}^m u_m^2dt$.
Similarly to find the state effort, we use the formula $\int_0^t  x^2dt$. As expected, since the proposed system uses four control inputs as opposed to three control inputs for the system in \cite{labar2019newton}, the effort
exerted by the controls for the proposed system is larger, however, the system as a whole is more efficient and exerts less energy. That is, the trade-off between faster convergence and larger control effort is not
balanced and shows the advantage of faster convergence in the system effort as a whole. However, for very large $t$ it is expected that the difference between both systems in state effort will diminish since the
proposed system is stabilized with a larger radius about the extremum.}
This drawback can be reduced/minimized if one considers techniques of attenuating oscillation \cite{pokhrel2023SIAMcontrol} especially given how fast the convergence is for the proposed ESC. Another option is the design
of a system whose control vanishes at the extremum point similar to what was proposed in \cite{VectorFieldGRUSHKOVSKAYA2018} which we tackle in Example 4.

\begin{figure}
\includegraphics{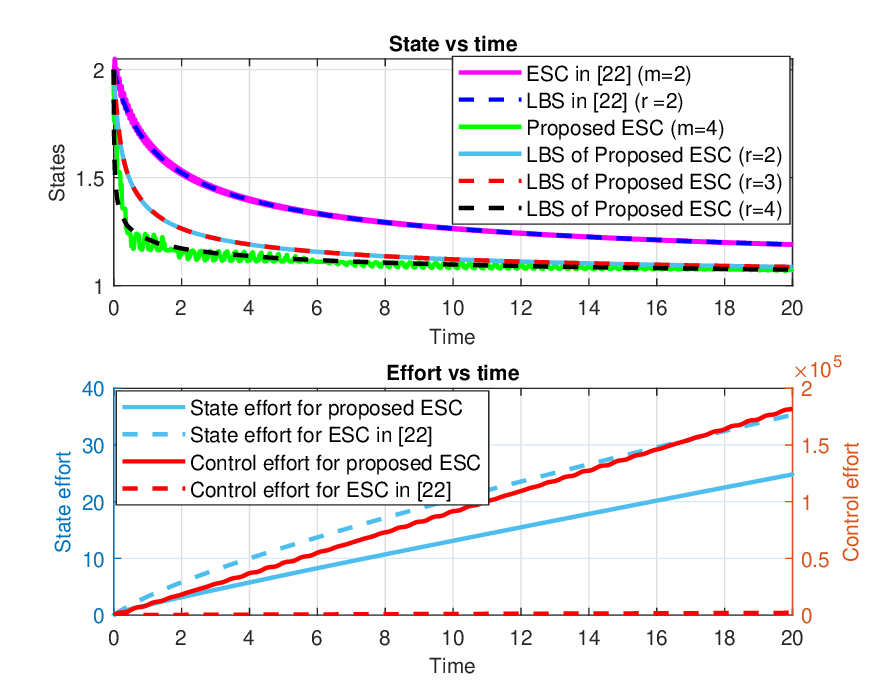}
\caption{\label{fig:vanishingClass}Comparison of performance between ESC with $m=2$ from \cite{VectorFieldGRUSHKOVSKAYA2018} and our proposed ESC with $m=4$ based on third-order LBS design.}
\end{figure}

\textbf{Example 4.} For this example, we take a system similar to the one in \cite[Section 4]{VectorFieldGRUSHKOVSKAYA2018}
which possesses {polynomial} convergence properties. {This system is a special case of the proposed system below when only the first two control inputs are considered with $p_1=p_2=0.5$:}
\begin{equation}
\label{eqn:vanishin_prop_main}
\begin{split}
\dot{x} = \sqrt{\frac{1-e^{-J(x)}}{1+e^{J(x)}}} \Bigm (\omega ^{p_1}sin(\psi) u_1 + \omega ^{p_2} cos(\psi)u_2  +\\ \omega ^{p_3} sin(\psi)u_3 + \omega ^{p_4} cos(\psi)u_4 \Bigm),
\end{split}
\end{equation}
for $J(x) \ne 0$ and $\dot{x} = 0$ for $J(x)=0$ with $\psi = e^{J(x)}+2 ln(e^{J(x)}-1)$ with $u_1 = cos(k_1\omega t)$, $u_2 = sin(k_2\omega t)$, $u_3 = cos(k_3\omega t)$, $u_4 = sin(k_4\omega t)$. We take $J = H(x-1)^4$. {We choose $p_1=p_3 = 0.99$, $p_2=p_4 = 0.01$, $k_1=k_2=1, k_3 = k_4 = 0.25$.} {We note that
for this example, similar to \cite{VectorFieldGRUSHKOVSKAYA2018}, Assumption A1 may be relaxed with a regularity requirement for the Lie derivatives which we assume is holding.} For simulation, we take $\omega = 100$,
$H=1/3$ and $x_0 = 2$. We use \eqref{eqn:LBS_main} to obtain a third-order LBS for the proposed ESC \eqref{eqn:vanishin_prop_main} and borrow the first-order LBS for the ESC with $m=2$ from \cite[Section 4]{VectorFieldGRUSHKOVSKAYA2018} with  $p_1=p_2=0.5$. {Here we briefly explain our design rationale, but the reader is directed to our supplementary file and code \cite{github} for all details. It is clear from
\eqref{eqn:vanishin_prop_main} that $b_1 = b_3$ and $b_2 = b_4$ with $b_1 = \sqrt{\phi^{*}}\sin{\psi}, b_2 = \sqrt{\phi^{*}} \cos{\psi}$, $\phi^{*} = (1-e^{-J(x)})/(1+e^{J(x)})$.
Now, for $L_2$ corresponding to first-order Lie brackets (see \eqref{eqn:l2_main}), we have $[b_1,b_3]=[b_2,b_4]=0$. So, $[b_1,b_2],[b_1,b_4],[b_2,b_3],[b_3,b_4]$ are the relevant non-zero Lie brackets (all of which depend on $[b_1,b_2]$); and only $\nu_{12},\nu_{14},\nu_{23},\nu_{34}$ are relevant. It was found that only $\nu_{12}$ and $\nu_{34}$ are non-zeros given by $\nu_{12}=0.5$ and $\nu_{34}=2$, both of which have $p_{j_1}+p_{j_2}=1$, yielding $L_2$ bounded as $\omega \rightarrow \infty$.
Now, for $L_3$ corresponding to second-order Lie brackets (see \eqref{eqn:l3_main}), Lie brackets having $[b_1,b_3]$ or $[b_2,b_4]$ in them will be zero and the relevant Lie brackets are $[b_1,[b_1,b_2]]$ and $[b_2,[b_1,b_2]]$. There are 16 relevant $\nu_{j_1j_2j_3}$ (corresponding to non-zero Lie brackets) given by $\bm{\nu}_{\bm{121}}$, $\nu_{122}$, $\bm{\nu}_{\bm{123}}$, $\nu_{124}$, 
        $\bm{\nu}_{\bm{141}}$, $\nu_{142}$, $\bm{\nu}_{\bm{143}}$, $\nu_{144}$,
        $\bm{\nu}_{\bm{231}}$, $\nu_{232}$, $\bm{\nu}_{\bm{233}}$, $\nu_{234}$, 
        $\bm{\nu}_{\bm{341}}$, $\nu_{342}$, $\bm{\nu}_{\bm{343}}$, $\nu_{344}$,
where bold $\nu_{j_1j_2j_3}$ correspond to those with $p_{j_1}+p_{j_2}+p_{j_3} \approx 2$. However, the iterated integrals for said $\nu_{j_1j_2j_3}$ were found to be zero. Thus, only $\nu_{j_1j_2j_3}$ with $p_{j_1}+p_{j_2}+p_{j_3} \approx 1$ (i.e.~with power of $\omega$ approximately about $-1$ in  \eqref{eqn:nuijk_main}) are the only ones with non-zero values, thus making $L_3$ bounded as $\omega \rightarrow \infty$, but less significant (i.e.,~$L_3$ almost vanish for large $\omega$). Nevertheless, none of the terms are ignored in the simulation.
For $L_4$ corresponding to third-order Lie brackets and extracting third-order derivative information (see \eqref{eqn:l4_main}) one will find that the Lie brackets of the form $[[b_{j_1},b_{j_2}],[b_{j_3},b_{j_4}]]$ are zeros and the non-zero relevant Lie brackets are $[[[b_1,b_2],b_1],b_1],[[[b_1,b_2],b_1],b_2],[[[b_1,b_2],b_2],b_1]$ and $[[[b_1,b_2],b_2],b_2]$. Since Lie brackets multiplied with $\beta_{1_{j_1j_2j_3j_4}}$ (see \eqref{eqn:l4_main}) are all zeros, $\beta_{1_{j_1j_2j_3j_4}}$ do not need to be computed and there are only 64 relevant $\beta_{2_{j_1j_2j_3j_4}}$ that need to be computed. All of the $\beta_{2_{j_1j_2j_3j_4}}$ have not been provided due to space constraints but it is found that $\beta_{2_{j_1j_2j_3j_4}}$ corresponding to those with $p_{j_1}+p_{j_2}+p_{j_3}+p_{j_4} \approx 3$ do not all have zero iterated integrals. Thus, $L_4$ is bounded as $\omega \rightarrow \infty$ and is significant in the higher-order LBS.} 
 \xref{fig:vanishingClass} {(top) shows a clear faster convergence advantage for the proposed ESC compared to the literature version \cite{VectorFieldGRUSHKOVSKAYA2018}. The proposed ESC is qualitatively captured by the
third-order LBS better than first- and second-order LBS. We present control effort and system state effort in  \xref{fig:vanishingClass} (bottom) using the same formulas as in Example 3. As expected, the proposed ESC
requires a much greater control effort as it has four control inputs instead of two control inputs in \cite{VectorFieldGRUSHKOVSKAYA2018}. However, the proposed system is much more efficient for the whole state effort
with accumulating advantage over time (unlike Example 3).   
  }

\section{Conclusion}

This paper provided a general and unifying theory that settles, connects, and bridges the concepts of higher-order Lie bracket approximation with higher-order averaging for a general class of control-affine systems \eqref{eqn:controlAffineIntro} {that include drift vector fields and not limited to underactuation conditions.} We conclude the equivalence between ($n$)-order LBS approximation and ($n+1$)-order averaging providing for the first time in literature a sequential framework for deriving higher-orders LBS approximating \eqref{eqn:controlAffineIntro}. Our results generalize, and indeed reduce to, previous efforts in the literature which were specific to less generalized classes of \eqref{eqn:controlAffineIntro} and concerned with particular orders of approximation. Consequently, the results of this paper can benefit the analysis/design of control-affine systems, including ESCs. For that, we provided a clearer understanding of why higher-order LBS (now higher-order averaging) may be needed and how higher-order LBS can become a complete averaging asymptote. Said concept provides insights for design improvements (e.g.,~much better convergence rate) of ESC and enables access for the first time in literature to third-order-derivative-based ESC designs. We provided multiple Examples and numerical simulations to demonstrate the effectiveness of our results. We hope this paper can bridge, benefit, and enable research between different communities of control theory and optimization.

\printauthorcredits

\appendix{}

\section[pfx={Appendix\space}]{Appendix: Proof of \xref{thm:mainthm_main}
} \label{sec:app}

  The term $I$ in \eqref{eqn:averaging_form_main} is independent of $\tau$, and its averaging is just itself. So, it can be later added directly to the averaged system. Only the term $II$ of \eqref{eqn:averaging_form_main} is considered. For simplicity and space considerations, the argument of $\bm{b}(\bm{x})$ is dropped, the integration variable $d \tau$ is dropped in $\int_0^{T'} (.,\tau)d\tau$, $\partial$ is used to denote $\partial/\partial \bm{x}$, and subscripts $i,j,e,l$ are used instead of $j_1,j_2,j_3,j_4$. We use $u'_i(k_i\tau)= \eta_i u_i(k_i \tau)$, and the notion $U'_i(k_i \tau) = \int_0^\tau u'_i(k_ip)dp $. With  \xref{lemma:wellposed_main} in place, \eqref{eqn:averaging_form_main} can be written as:
\numberwithin{equation}{section}\setcounter{equation}{0}
\begin{equation}
\label{eqn:dxdtau1_app}
    \frac{d\bm{x}}{d\tau} = \epsilon \sum_{i=1}^m u'_i(k_i\tau)\bm{b}_i(\bm{x})\text{},
\end{equation}
with $\epsilon = \omega^{p^{*}-1}$. Next, in $\tau$- scale we find $\stmboldsymbol{\bLambda}_r$ corresponding to $r=1,2,3,\text{ and }4$ given by \eqref{eqn:lambda1_main}--\eqref{eqn:lambda4_main}. Then, we convert $\stmboldsymbol{\bLambda}_r$ to $t$- scale which we denote by $\bm{L}_r$.

\textbf{For} $\bm{r = 1}$, we have
\begin{equation}
\stmboldsymbol{\bLambda}_1 = \frac{\epsilon}{T'}\int_0^{T'}  \sum_{i=1}^m u'_i(k_i\tau)\bm{b}_i(\bm{x}).
\end{equation}
From A2, $u'_i(k_i\tau)$ is $T'-$ periodic with zero average. It follows that $\stmboldsymbol{\bLambda}_1=0$. Hence, the first-order averaging term in \eqref{eqn:averaged_prem} vanishes (i.e.,~$ \stmboldsymbol{\bLambda}_1 = 0$). Similarly, the averaged term in $t$- scale is $\bm{L}_1 = \omega  \stmboldsymbol{\bLambda}_1 = 0$.

\textbf{For} $\bm{r = 2}$, let us consider
\begin{align}
\bm{A} &= \bm{f}_\tau = \sum_{i=1}^m u'_i(k_i \tau) \bm{b}_i,\label{eqn:A_app}\\
    \bm{B} &= \int_0^\tau \bm{f}_p dp = \sum_{i=1}^mU'_i(k_i \tau) \bm{b}_i \label{eqn:B_app}\\
    \stmboldsymbol{\bLambda}'_2 &= \frac{1}{T'}\int_0^{T'} [\bm{B},\bm{A}]\label{eqn:lambda2_1_app}
\end{align}
     where obviously $[\bm{B},\bm{A}]={\partial \bm{A}} \; \bm{B}-\partial \bm{B} \; \bm{A}$. Furthermore,  
         \begin{equation*}
\begin{split}
&[\bm{B},\bm{A}] = \sum_{i = 1}^m u'_{i}( k_{i}\tau)\partial \bm{b}_{i}
    \sum_{j = 1}^mU'_{j}(k_{j} \tau) \bm{b}_{j}
\end{split}
\end{equation*}
\begin{equation}
\label{eqn:ba_app}
\begin{split}
&-\sum_{i = 1}^m U'_{i}(k_{i} \tau) \partial \bm{b}_{i}  \sum_{j = 1}^m u'_{j}(k_{j}\tau) \bm{b}_{j}\\
    &= \sum_{i = 1}^m \sum_{j = 1}^m \Bigl(u'_{i}(k_{i}\tau)U'_{j}(k_{j}\tau)\partial \bm{b}_{i}\bm{b}_{j}-
     u'_{j}(k_{j}\tau)U'_{i}(k_{i}\tau)\partial \bm{b}_{i}\bm{b}_{j} \Bigr)\\
    &= \sum_{i = 1}^m \sum_{j=i+1}^m \Bigl( u'_{j}(k_{j}\tau)U'_{i}(k_{i}\tau) -
    u'_{i}(k_{i}\tau)U'_{j}(k_{j}\tau) \Bigr) [\bm{b}_{i}, \bm{b}_{j}].
\end{split}
\end{equation}
Now, from \eqref{eqn:lambda2_1_app} and \eqref{eqn:ba_app}, we have:
\begin{equation}
\begin{split}
\stmboldsymbol{\bLambda}'_2 &= \frac{1}{T'}\int_0^{T'}[\bm{B},\bm{A}]
    =  \sum_{i=1}^m \sum_{j=i+1}^m \nu'_{ij}[\bm{b}_{i}, \bm{b}_{j}],
\end{split}
\end{equation}
with $\nu'_{i j} = \frac{1}{T'}\int_0^{T'}  \left( u'_{j}(k_{j}\tau)U'_{i}(k_{i}\tau) -u'_{i}(k_{i}\tau)U'_{j}(k_{j}\tau) \right)$.
Now, the   second-order averaging term in \eqref{eqn:averaged_prem} is:
\begin{equation}
\begin{split}
\stmboldsymbol{\bLambda}_2 = \frac{\epsilon^2}{2} \stmboldsymbol{\bLambda}'_2 = \frac{\epsilon^2}{2} \sum_{i=1}^m \sum_{j=i+1}^m \nu'_{i j}[\bm{b}_{i}, \bm{b}_{j}].
\end{split}
\end{equation}
Now, going back to $t$-scale and substituting $\epsilon = \omega^{p^{*}-1}$:
\begin{equation}
\label{eqn:nu_ij_dash_app}
\bm{L}_2= \omega \frac{\epsilon^2}{2} \stmboldsymbol{\bLambda}'_2 = \frac{\omega \omega^{2p^{*}-2}}{2} \sum_{i=1}^m \sum_{j=i+1}^m \nu'_{i j}[\bm{b}_{i}, \bm{b}_{j}],
\end{equation}
Now, using the relations $u'_i(k_i \tau) = \eta_i u_i(k_i \tau)$  and $\eta_i = \omega^{p_i-p^{*}}$ from  \xref{lemma:wellposed_main}, \eqref{eqn:nu_ij_dash_app} can be written as
\begin{equation}
\bm{L}_2 =  \sum_{i=1}^m \sum_{j=i+1}^m \nu_{i j}[\bm{b}_{i}, \bm{b}_{j}],
\end{equation}
with
\begin{equation}
\nu_{i j} = \frac{\omega^{p_{i}+p_{j}-1}}{2T'} \int_0^{T'}  \left( u_{j}(k_{j}\tau)U_{i}(k_{i}\tau) -u_{i}(k_{i}\tau)U_{j}(k_{j}\tau) \right).
\end{equation}
 \textbf{For} $\bm{r = 3}$, since $\stmboldsymbol{\bLambda}_1 = 0$, then $\stmboldsymbol{\bLambda}_3 $ in \eqref{eqn:lambda3_main} is written as
\begin{equation}
\label{eqn:lambda3_app}
    \stmboldsymbol{\bLambda}_3 = \frac{\epsilon^3}{3!} \stmboldsymbol{\bLambda}_3'= \frac{2\epsilon^3}{3!T'}\int_0^{T'}[\bm{B},[\bm{B},\bm{A}]].
\end{equation}
So, we first need to find $[\bm{B},[\bm{B},\bm{A}]]$, but before that we introduce the following notation for the sake of simplicity
\begin{equation}
\beta'_{ij} =  \left( u'_j(k_j\tau)U'_i(k_i\tau) -u'_i(k_i\tau)U'_j(k_j\tau) \right).
\end{equation}
Then, $[\bm{B},[\bm{B},\bm{A}]]$ can be found as
\begin{equation*}
\begin{split}
&[\bm{B},[\bm{B},\bm{A}]] = \partial[\bm{B},\bm{A}]\bm{B} - \partial \bm{B}[\bm{B},\bm{A}]\\
        &=  \sum_{i=1}^m \sum_{j=i+1}^m \beta'_{i j} \partial [\bm{b}_{i},\bm{b}_{j}] \sum_{e=1}^m \bm{b}_{e} U'_{e}(k_{e}\tau)
\end{split}
\end{equation*}
\xpar\cvskip[-1pc]
\begin{equation}
\label{eqn:bba_app}
\begin{split}
&- \sum_{i=1}^m\partial \bm{b}_{i}U'_{i} (k_{i}\tau) \sum_{j = 1}^m \sum_{e= j+1}^m \beta'_{j e} [\bm{b}_{j},\bm{b}_{e}]\\
        & =  \sum_{i=1}^m \sum_{j=i+1}^m \sum_{e=1}^m \beta'_{i j} U'_{e}(k_{e}\tau) \partial [\bm{b}_{i},\bm{b}_{j}] \bm{b}_{e}\\
&-\xs  \sum_{i=1}^m  \sum_{j = 1}^m \sum_{e= j+1}^m U'_{i}(k_{i}\tau)\beta'_{j e}\partial \bm{b}_{i}[\bm{b}_{j},\bm{b}_{e}]\\
        &=  \sum_{i=1}^m \sum_{j=i+1}^m \sum_{e=1}^m U'_{e}(k_{e}\tau) \beta'_{i j} [\bm{b}_{e},[\bm{b}_{i},\bm{b}_{j}]]
\end{split}
\end{equation}

Now, using \eqref{eqn:bba_app} and \eqref{eqn:lambda3_app}, we get
\begin{equation}
\begin{split}
&\stmboldsymbol{\bLambda}'_3   =  \frac{1}{T'}\int_0^{T'}  \sum_{i=1}^m \sum_{j=i+1}^m \sum_{e=1}^m \beta'_{i j}U'_{e} (k_{e} \tau) [\bm{b}_{e},[\bm{b}_{i},\bm{b}_{j}]]\\
    & = \sum_{i=1}^m \sum_{j=i+1}^m \sum_{e=1}^m \nu'_{i j e} [\bm{b}_{e},[\bm{b}_{i},\bm{b}_{j}]]
\end{split}
\end{equation}
with $\nu'_{i j e} = \frac{1}{T'}\int_0^{T'}  \beta'_{i j}U'_e(k_e \tau) d\tau$. Thus, the term corresponding to third-order averaging in \eqref{eqn:averaged_prem} is
\begin{equation}
\begin{split}
{\stmboldsymbol{\bLambda}}_{\bm{3}} = \frac{\epsilon^3 \stmboldsymbol{\bLambda}'_3}{3}
        & = \frac{\epsilon^3}{3}\sum_{i=1}^m \sum_{j=i+1}^m \sum_{e=1}^m \nu'_{i j e} [\bm{b}_{e},[\bm{b}_{i},\bm{b}_{j}]].
\end{split}
\end{equation}
Now, in $t$- scale with $\epsilon = \omega^{p^{*}-1}$ from  \xref{lemma:wellposed_main}, we get
\begin{equation}
\label{eqn:L3_dash_app}
\begin{split}
\bm{L}_3 &= \omega \frac{\epsilon^3 \stmboldsymbol{\bLambda}_3}{3}\\
    &=  \frac{1}{3} \omega \omega^ {3p^{*}-3}\sum_{i=1}^m \sum_{j=i+1}^m \sum_{e=1}^m \nu'_{i j e} [\bm{b}_{e},[\bm{b}_{i},\bm{b}_{j}]],
\end{split}
\end{equation}
Now, using the relations $u'_i(k_i \tau) = \eta_i u_i(k_i \tau)$  and $\eta_i = \omega^{p_{i}-p^{*}}$ from  \xref{lemma:wellposed_main}, \eqref{eqn:L3_dash_app} can be written as
\begin{equation}
\bm{L}_3 =  \sum_{i=1}^m \sum_{j=i+1}^m \sum_{e=1}^m  \nu_{i j e} [\bm{b}_{e},[\bm{b}_{i},\bm{b}_{j}]]
\end{equation}
with
\begin{align*}
\nu_{i j e} &= \frac{\omega^{p_{i}+p_{j}+p_{e}-2}}{3T'}\int_0^{T'} \biggm((u_{j}(k_{j}\tau)U_{i}(k_{i}\tau) \\
&-\xs  u_{i}(k_{i}\tau)U_{j}(k_{j}\tau) ) U_{e} (k_{e} \tau) \biggm)d\tau.
\end{align*}

\textbf{For} $\bm{r=4}$, we introduce the following notations. Symbols $\bm{A}$ and $\bm{B}$ are defined as in \eqref{eqn:A_app}--\eqref{eqn:B_app} but with a subscript that denotes the argument of $u'$ and $U'$. Now, we aim at writing $\stmboldsymbol{\bLambda}_4$ from \eqref{eqn:lambda4_main} as $\stmboldsymbol{\bLambda}_4 = 2\epsilon^4\stmboldsymbol{\bLambda}'_4/4!$ with
\begin{equation}
\label{eqn:lambda4_app}
\begin{split}
\stmboldsymbol{\bLambda}'_4    = \frac{1}{T'}\int_0^{T'} \biggl( \int_0^\tau [\bm{C},\bm{D}] dp +[\bm{E},\bm{f}_\tau] + \int_0^\tau \bm{G} dp \biggr).
\end{split}
\end{equation}
In the next step, we find the expressions for $\bm{C},\bm{D},\bm{E},\bm{G}$ and associated brackets.
Now, for the computation of $\bm{C},\bm{D}$ and the first part of \eqref{eqn:lambda4_app}, we have
\begin{align*}
\nf&[\bm{B}_q,\bm{A}_q] = \sum_{i=1}^m \partial \bm{b}_{i} u'_{i}(k_i q) \sum_{j=1}^m \bm{b}_{j} U'_{j}(k_j q) \\ 
&-  \sum_{i=1}^m \partial \bm{b}_{i} U'_{i}(k_i q) \sum_{j=1}^m \bm{b}_{j} u'_{j}(k_j q) \\
    &= \sum_{i=1}^m \sum_{j=i+1}^m (u'_{j}(k_j q)U'_{i}(k_i q)-u'_{i}(k_i q)U'_{j}(k_j q))[\bm{b}_{i},\bm{b}_{j}]\\
    &= \sum_{i=1}^m \sum_{j=i+1}^m \alpha'_1(q)_{i j} [\bm{b}_{i},\bm{b}_{j}].
\end{align*}
Now, we can compute $\bm{C}$ as
\begin{equation}
\bm{C} = \int_0^p [\bm{B}_q,\bm{A}_q] = \sum_{i=1}^m \sum_{j=i+1}^m \alpha'_2(p)_{i j} [\bm{b}_{i},\bm{b}_{j}],
\end{equation}
where $\alpha'_2(p)_{i j}= \int_0^p \alpha'_1(q)_{i j} dq$. Next, we compute $\bm{D}$ as
\begin{equation}
\begin{split}
&\bm{D} = \left[\sum_{i=1}^m \bm{b}_{i} u'_{i}(p),\sum_{i=1}^m \bm{b}_{i} u'_{i} (\tau) \right]\\
    &= \sum_{i=1}^m u'_{i}(\tau) \partial \bm{b}_{i} \sum_{j=1}^m u'_{j}(p) \bm{b}_{j} - \sum_{i=1}^m u'_{i}(p) \partial \bm{b}_{i} \sum_{j=1}^m u'_{j} (\tau) \bm{b}_{j}\\
    & = \sum_{i=1}^m \sum _{j=i+1}^m (u'_{j}(\tau) u'_{i}(p) - u'_{i}(\tau) u'_{j}(p))[\bm{b}_{i}, \bm{b}_{j}]\\
    & = \sum_{i=1}^m \sum _{j=i+1}^m \alpha'_3(\tau,p)_{i j} [\bm{b}_{i}, \bm{b}_{j}].
\end{split}
\end{equation}
Next, the Lie bracket of $\bm{C}$ and $\bm{D}$ can be found as
\begin{equation*}
\begin{split}
&[\bm{C},\bm{D}] \\
    &= \left[\sum_{i=1}^m \sum_{j=i+1}^m \alpha'_2(p)_{i j} [\bm{b}_{i},\bm{b}_{j}], \sum_{i=1}^m \sum _{j=i+1}^m \alpha'_3(\tau,p)_{i j} [\bm{b}_{i}, \bm{b}_{j}] \right]\\
    &= \sum_{i=1}^m \sum_{j=i+1}^m \alpha'_3(\tau,p)_{i j} \partial [\bm{b}_{i},\bm{b}_{j}] \sum_{e=1}^m \sum_{l=e+1}^m \alpha'_2(p)_{e l} [\bm{b}_{e},\bm{b}_{l}]\\
&-\xs  \sum_{i=1}^m \sum_{j=i+1}^m \alpha'_2(p)_{i j} \partial [\bm{b}_{i},\bm{b}_{j}] \sum_{e=1}^m \sum_{l=e+1}^m \alpha'_3(\tau,p)_{e l} [\bm{b}_{e},\bm{b}_{l}]\\
    &= \sum_{i=1}^m \sum_{j=i+1}^m \sum_{e=1}^m \sum_{l=e+1}^m \alpha'_4(\tau,p)_{i j e l} \Bigl[[\bm{b}_{i},\bm{b}_{j}],[\bm{b}_{e},\bm{b}_{l}]\Bigr]
\end{split}
\end{equation*}
where $\alpha'_4(\tau,p)_{i j e l} = \alpha'_3(\tau,p)_{e l} \alpha'_2(p)_{i j}$ and $(i,j)\ne (e,l)$.
Finally, we integrate $[\bm{C},\bm{D}]$ from 0 to $\tau$ as
\begin{equation}
\begin{split}
\int_0^\tau [\bm{C},\bm{D}] =& \sum_{i=1}^m \sum_{j=i+1}^m \sum_{e=1}^m \sum_{l=e+1}^m \alpha'_5(\tau)_{i j e l} \\
    &\Bigl[[\bm{b}_{i},\bm{b}_{j}],[\bm{b}_{e},\bm{b}_{l}]\Bigr]
\end{split}
\end{equation}
with $\alpha'_5(\tau)_{i j e l} = \int_0^\tau \alpha'_4(\tau,p)_{i j e l} dp$.

Now, for the computation of second part of \eqref{eqn:lambda4_app}, first we find
\begin{align*}
\nf&[\bm{C},\bm{f}_p] = \biggl[ \sum_{i=1}^m \sum_{j=i+1}^m \alpha'_2 (p)_{i j} [\bm{b}_{i},\bm{b}_{j}], \sum_{i=1}^m \bm{b}_{i} u'_{i}(p) \biggr]\\
    & = \sum_{i=1}^m \partial \bm{b}_{i} u'_{i}(p) \sum_{j=1}^m \sum_{e=j+1}^m \alpha'_2(p)_{j e}  [\bm{b}_{j},\bm{b}_{e}] \\
&-\xs  \sum_{i=1}^m \sum_{j=i+1}^m \alpha'_2(p)_{i j} \partial [\bm{b}_{i},\bm{b}_{j}]\sum_{e=1}^m  \bm{b}_{e} u'_{e}(p),\\
    &= \sum_{i=1}^m \sum_{j=1}^m \sum_{e=j+1}^m  u'_{i}(p)\alpha'_2(p)_{j e} \partial \bm{b}_{i} [\bm{b}_{j},\bm{b}_{e}]\\
&-\xs  \sum_{i=1}^m \sum_{j=i+1}^m \sum_{e=1}^m 
    \alpha'_2(p)_{i j} u'_{e}(p) \partial [\bm{b}_{i},\bm{b}_{j}]\bm{b}_{e} \\
    &= \sum_{i=1}^m \sum_{j=i+1}^m \sum_{e=1}^m \alpha'_2(p)_{i j} u'_{e}(p) \biggl[[\bm{b}_{i},\bm{b}_{j}],\bm{b}_{e}\biggr],
\end{align*}
Thus,
\begin{equation*}
\begin{split}
[\bm{C},\bm{f}_p] = \sum_{i=1}^m \sum_{j=i+1}^m \sum_{e=1}^m \alpha'_6(p)_{i j e} \biggl[[\bm{b}_{i},\bm{b}_{j}],\bm{b}_{e}\biggr],
\end{split}
\end{equation*}
with $\alpha'_6(p)_{i j e} = \alpha'_2(p)_{i j} u'_{e}(p)$, $\bm{E}$ can be computed as
\begin{equation}
\bm{E} = \int_0^\tau [\bm{C},\bm{f}_p] = \sum_{i=1}^m \sum_{j=i+1}^m \sum_{e=1}^m \alpha'_7(\tau)_{i j e} \biggl[[\bm{b}_{i},\bm{b}_{j}],\bm{b}_{e}\biggr],
\end{equation}
with  $\alpha'_7(\tau)_{i j e} = \int_0^\tau \alpha'_6(p)_{i j e} dp$. Additionally,
\begin{equation*}
\begin{split}
[\bm{E},\bm{f}_\tau] &= \Biggl[ \sum_{i=1}^m \sum_{j=i+1}^m \sum_{e=1}^m \alpha'_7(\tau)_{i j e}[[\bm{b}_{i},\bm{b}_{j}],\bm{b}_{e}] \sum_{i=1}^m u'_{i}(\tau) \bm{b}_{i} \Biggr],\\
    &= \sum_{i=1}^m u'_{i}(\tau) \partial \bm{b}_{i} \sum_{j=1}^m \sum_{e=j+1}^m \sum_{l=1}^m \alpha'_7(\tau)_{j e l}[[\bm{b}_{j},\bm{b}_{e}],\bm{b}_{l}] \\
&-\xs  \sum_{i=1}^m \sum_{j=i+1}^m \sum_{e=1}^m \alpha'_7(\tau)_{i j e} \partial [[\bm{b}_{i},\bm{b}_{j}]\bm{b}_{e}] \sum_{l=1}^m u'_{l} (\tau)\bm{b}_{l},
\end{split}
\end{equation*}
\begin{equation}
\begin{split}
& = \sum_{i=1}^m \sum_{j=i+1}^m \sum_{e=1}^m \sum_{l=1}^m \alpha'_7(\tau)_{i j e}u'_{l} (\tau)[[[\bm{b}_{i},\bm{b}_{j}],\bm{b}_{e}],\bm{b}_{l}],\\
    & = \sum_{i=1}^m \sum_{j=i+1}^m \sum_{e=1}^m \sum_{l=1}^m \alpha'_8(\tau)_{i j e l} [[[\bm{b}_{i},\bm{b}_{j}],\bm{b}_{e}],\bm{b}_{l}],
\end{split}
\end{equation}
where $\alpha'_8(\tau)_{i j e l} = \alpha'_7(\tau)_{i j e} u'_{l}(\tau)$.
Finally, for the computation of the third part of \eqref{eqn:lambda4_app}, we define $\bm{F}$ as
\begin{equation}
\begin{split}
&\bm{F}  = \left[\left[ \int_0^p \bm{f}_q dq, \bm{f}_p\right],\bm{A}_\tau\right]=  \left[\left[ \bm{B}_p, \bm{A}_p\right],\bm{A}_\tau\right]\\
    & = \left[ \sum_{i=1}^m \sum_{j=i+1}^m \alpha'_1(p)_{i j} [\bm{b}_{i},\bm{b}_{j}], \sum_{i=1}^m \bm{b}_{i} u'_{i}(\tau) \right]\\
    &= \sum_{i=1}^m \sum_{j=i+1}^m \sum_{e=1}^m \alpha'_1(p)_{i j} u'_{e}(\tau) [[\bm{b}_{i},\bm{b}_{j}],\bm{b}_{e}]\\
    &= \sum_{i=1}^m \sum_{j=i+1}^m \sum_{e=1}^m \alpha'_9(p,\tau)_{i j e} [[\bm{b}_{i},\bm{b}_{j}],\bm{b}_{e}],
\end{split}
\end{equation}
with $\alpha'_9(p,\tau)_{i j e} = \alpha'_1(p)_{i j} u'_{e}(\tau)$. Next, we find $\bm{G}$ as
\begin{equation}
\begin{split}
\bm{G} &= \biggl[\sum_{i=1}^m \bm{b}_{i} U'_{i} (p), \bm{F} \biggr]\\
    & = \biggl[\sum_{i=1}^m \bm{b}_{i} U'_{i}(p), \sum_{i=1}^m \sum_{j=i+1}^m \sum_{e=1}^m \alpha'_9 (p,\tau)_{i j e} [[\bm{b}_{i},\bm{b}_{j}], \bm{b}_{e}] \biggr],\\
    & = \sum_{i=1}^m \sum_{j=i+1}^m \sum_{e=1}^m\alpha_9'(p,\tau)_{i j e} \partial[[\bm{b}_{i},\bm{b}_{j}],\bm{b}_{e}] \sum_{l=1}^m \bm{b}_{l} U'_{l} (p)\\
&-\xs  \sum_{i=1}^m \partial \bm{b}_{i} U'_{i}(p) \sum_{j=1}^m \sum_{e= j+1}^m \sum_{l=1}^m [[\bm{b}_{j},\bm{b}_{e}], \bm{b}_{l}] \alpha_9'(p,\tau)_{j e l},\\
    & = \sum_{i=1}^m \sum_{j=i+1}^m \sum_{e=1}^m \sum_{l=1}^m \alpha_9' (p,\tau)_{i j e}U'_{l}(p) [\bm{b}_{l},[[\bm{b}_{i},\bm{b}_{j}],\bm{b}_{e}]].
\end{split}
\end{equation}
We can now find the integral of $\bm{G}$ from 0 to $\tau$ as
\begin{equation*}
    \int_0^\tau \bm{G}=  \sum_{i=1}^m \sum_{j=i+1}^m \sum_{e=1}^m \sum_{l=1}^m \alpha'_{10}(\tau)_{i j e l} [\bm{b}_{l},[[\bm{b}_{i},\bm{b}_{j}],\bm{b}_{e}]],
\end{equation*}
where $\alpha'_{10}(\tau)_{i j e l} =  \int_0^\tau \alpha'_9(p,\tau)_{i j e} U'_{l} (p) dp $.
Thus, $\stmboldsymbol{\bLambda}'_4$ from \eqref{eqn:lambda4_app}, after computing $\bm{C},\bm{D},\bm{E},\bm{G}$ and associated brackets, can be written as
\begin{equation*}
\begin{split}
\stmboldsymbol{\bLambda}'_4 &= \frac{1}{T'}\int_0^{T'} \Biggm(\sum_{i=1}^m \sum_{j=i+1}^m \sum_{e=1}^m \sum_{l=e+1}^m \alpha'_5(\tau)_{i j e l} 
    \Bigl[[\bm{b}_{i},\bm{b}_{j}],[\bm{b}_{e},\bm{b}_{l}]\Bigr]
\end{split}
\end{equation*}
\begin{equation}
\begin{split}
&+ \sum_{i=1}^m \sum_{j=i+1}^m \sum_{e=1}^m \sum_{l=1}^m \alpha'_8(\tau)_{i j e l} [[[\bm{b}_{i},\bm{b}_{j}],\bm{b}_{e}],\bm{b}_{l}]\\
&+\xs  \sum_{i=1}^m \sum_{j=i+1}^m \sum_{e=1}^m \sum_{l=1}^m \alpha'_{10}(\tau)_{i j e l} [\bm{b}_{l},[[\bm{b}_{i},\bm{b}_{j}],\bm{b}_{e}]] \Biggm).
\end{split}
\end{equation}
Finally, changing back to $t$-scale, and using $\epsilon = \omega^{p^{*}-1}$ we get
\begin{equation}
\begin{split}
& \bm{L}_4 =  \omega \frac{\omega^{4p^{*}-4} \stmboldsymbol{\bLambda}'_4}{12} = \frac{\omega^{4p^{*}-3} \stmboldsymbol{\bLambda}'_4}{12}\\
    & = \frac{\omega^{4p^{*}-3}}{12}\frac{1}{T'}\int_0^{T'} \Biggm(\sum_{i=1}^m \sum_{j=i+1}^m \sum_{e=1}^m \sum_{l=e+1}^m \alpha'_5(\tau)_{i j e l} \Bigl[[\bm{b}_{i},\bm{b}_{j}],[\bm{b}_{e},\bm{b}_{l}]\Bigr]\\
&+\xs  \sum_{i=1}^m \sum_{j=i+1}^m \sum_{e=1}^m \sum_{l=1}^m \alpha'_8(\tau)_{i j e l} [[[\bm{b}_{i},\bm{b}_{j}],\bm{b}_{e}],\bm{b}_{l}]\\
&+\xs  \sum_{i=1}^m \sum_{j=i+1}^m \sum_{e=1}^m \sum_{l=1}^m \alpha'_{10}(\tau)_{i j e l} [\bm{b}_{l},[[\bm{b}_{i},\bm{b}_{j}],\bm{b}_{e}]] \Biggm),\\ 
    & = \sum_{i=1}^m \sum_{j=i+1}^m \sum_{e=1}^m \sum_{l=e+1}^m \beta_{1_{i j e l}} \Bigl[[\bm{b}_{i},\bm{b}_{j}],[\bm{b}_{e},\bm{b}_{l}]\Bigr]\\
&+\xs  \sum_{i=1}^m \sum_{j=i+1}^m \sum_{e=1}^m \sum_{l=1}^m (\beta_{2_{i j e l}} -\beta_{3_{i j e l}})[[[\bm{b}_{i},\bm{b}_{j}],\bm{b}_{e}],\bm{b}_{l}]
\end{split}
\end{equation}
where
\begin{align*}
\beta_{1_{i j e l}} &=  \omega^{4p^{*}-3} \frac{1}{12T'}\int_0^{T'} \eta_{i} \eta_{j} \eta_{e} \eta_{l} \alpha_5(\tau)_{i j e l}\\
    &= \omega^{4p^{*}-3} \omega^{p_{i}+p_{j}+p_{e}+p_{l}-4p^{*}} \frac{1}{12T'}\int_0^{T'}  \alpha_5(\tau)_{i j e l}\\
    &= \omega^{p_{i}+p_{j}+p_{e}+p_{l}-3} \frac{1}{12T'}\int_0^{T'}  \alpha_5(\tau)_{i j e l}.
\end{align*} 
Similarly,
\begin{equation*}
\beta_{2_{i j e l}}= \omega^{p_{i}+p_{j}+p_{e}+p_{l}-3} \frac{1}{12T'}\int_0^{T'} \alpha_8(\tau)_{i j e l},
\end{equation*}
\begin{equation*}
\beta_{3_{i j e l}}= \omega^{p_{i}+p_{j}+p_{e}+p_{l}-3}  \frac{1}{12T'}\int_0^{T'} \alpha_{10}(\tau)_{i j e l}.
\end{equation*} 
By using the variable $\bm{z}$ instead of $\bm{x}$, the fourth-order averaging system in $t$-scale including the drift term $\bm{b}_0(\bm{z})$ is finally written as:
\begin{equation}
\dot{\bm{z}} = \bm{b}_{\bm{0}}(\bm{z})+ \bm{L}_1(\bm{z})+ \bm{L}_2 (\bm{z})+ \bm{L}_3(\bm{z}) + \bm{L}_4(\bm{z}).
\end{equation}
This completes the proof of \xref{thm:mainthm_main}. \qed
\bibliographystyle{ieeetr}
\bibliography{references}           

\end{document}